\documentclass[10pt]{article}  
\usepackage{amsmath, amsfonts, amsthm, amssymb, amscd}
\usepackage[mathcal]{euscript}
\usepackage{dsfont, pxfonts, mathrsfs, accents, verbatim}
\theoremstyle{plain}
        \newtheorem{theorem}{Theorem}[section]
\newtheorem{definition}[theorem]{Definition}
        \newtheorem{proposition}[theorem]{Proposition}
        \newtheorem{lemma}[theorem]{Lemma}
        \newtheorem{corollary}[theorem]{Corollary}

\numberwithin{equation}{section}
\newcommand \FC {{\mathcal FC}}
\newcommand \R      {\mathbb R}
\newcommand \be     {\begin{equation}}
\newcommand \ee     {\end{equation}}
 
\newcommand \del        \partial
\newcommand \eps        \varepsilon
\newcommand \auth   \textsc

\newcommand \Lie    {{\mathcal L}}
\newcommand \Dcal   {{\mathcal D}}
\newcommand \Jcal   {{\mathcal J}}
\newcommand \Ncal   {{\mathcal N}}
\newcommand \Acal   {{\mathcal A}}
\newcommand \Hcal   {{\mathcal H}}
\newcommand \Bcal   {{\mathcal B}}
\newcommand \Ical   {{\mathcal I}}
\newcommand \gt     {{\widetilde g}} 
\newcommand \Inj        {\text {\bf Inj}}
\newcommand \vol     {\text {\bf Vol}}
\newcommand \Conj    {\text {\bf Conj}}

\newcommand \Riem       {\text{\bf Rm}} 
\newcommand \Ric       {\text{\bf Ric}}
\newcommand \Gammat  {\widetilde{\Gamma}} 
\newcommand \Rt         {\widetilde{R}}
\newcommand \sigmat         {\widetilde{\sigma}}
 
\newcommand \expb   {{{\text{\bf exp}}}} 

\newcommand \Vt     {{\widetilde V}}
\newcommand \NullInj   {\text{\bf Null Inj}}
\newcommand \NullConj   {\text{\bf Null Conj}}
\newcommand \tsigma     {\widetilde \sigma}
\newcommand \Cbf    {\mathbf C}
\newcommand \Lbf    {\mathbf L}

\begin{document}

\title{Injectivity Radius of Lorentzian Manifolds}
\author{Bing-Long Chen\footnote{
Department of Mathematics, Sun Yet Sen University, 510275 Guang-Zhou, P.R. of China. E-mail: {\sl mcscbl@mail.sysu.edu.cn.}
}
\, and
%
Philippe G. LeFloch\footnote{Laboratoire Jacques-Louis Lions \& Centre National de la Recherche Scientifique, 
Universit\'e de Paris 6, 4 place Jussieu, 
75252 Paris, France. E-mail : {\sl lefloch@ann.jussieu.fr}
\newline
2000\textit{\ AMS Subject Classification:} 53C50, 83C05, 53C12, 53C22.
}
}
\date{\today}
\maketitle

\begin{abstract}
Motivated by the application to spacetimes of general relativity we investigate the geometry and regularity of
Lorentzian manifolds under certain curvature and volume bounds. 
We establish several injectivity radius estimates at a point or on the past null cone of a point. 
Our estimates are entirely local and geometric, and 
are formulated via a reference Riemannian metric that we canonically associate
with a given observer $(p,T)$ --where $p$ is a point of the manifold
and $T$ is a future-oriented time-like unit vector prescribed at $p$. 
The proofs are based on a generalization of arguments from Riemannian geometry.
We first establish estimates on the reference Riemannian metric, 
and then express them in term of the Lorentzian metric. In the context of general relativity, our estimates 
should be useful to investigate the regularity of spacetimes satisfying Einstein field equations.
\end{abstract}


\section{Introduction}
\label{IN-0}

\subsection*{Aims of this paper}
The regularity and compactness of Riemannian manifolds under a~priori bounds
on geometric quantities such as curvature, volume, or diameter represent important issues in Riemannian geometry.
In particular, the derivation of lower bounds on the injectivity radius of a Riemannian manifold,
and the construction of local coordinate charts in which the metric has optimal regularity
are now well-understood. Moreover,
Cheeger-Gromov's theory provides geometric conditions for the strong compactness of sequences of manifolds 
and has become a central tool in Riemannian geometry. See for instance \cite{Anderson0,Besse,Cheeger1,CheegerGromov1,CheegerGromov2,CFG,JostKarcher,Peters,Petersen}. 

Our objective in the present paper is to present some extension of these classical techniques and
results to Lorentzian manifolds. Recall that a Lorentzian metric is not positive definite,
but has signature $(-, +, \ldots, +)$.
Motivated by recent work by Anderson \cite{Anderson1} and Klainerman and Rodnianski \cite{KR4},
we derive here several injectivity radius estimates for Lorentzian manifolds satisfying
certain curvature and volume bounds. That is, we provide sharp lower bounds on the size of the geodesic ball
around one point within which the exponential map is a global diffeomorphism and, therefore, 
we obtain sharp control of the manifold geometry.
Our proofs rely on arguments that are known to be flexible and efficient in Riemannian geometry,
and are here extended to the Lorentzian setting; we analyze the properties of 
Jacobi fields and rely on volume comparison and homotopy arguments.

In our presentation
(see for instance our main result in Theorem~\ref{nofoliation} at the end of this introduction)
we emphasize the importance of having assumptions and estimates that 
are stated locally and geometrically, and avoid direct use of coordinates.
When necessary, coordinates should be constructed a~posteriori, 
once uniform bounds on the injectivity radius have been established.

Our motivation comes from general relativity, where one of the most challenging problems is the formation
and the structure of singularities in solutions to the Einstein field equations.
Relating curvature and volume bounds to the regularity of the manifold, as we do in this paper,
is necessary before tackling an investigation of the geometric properties of singular spacetimes
satisfying Einstein equations.
(See, for instance, \cite{Anderson1,Anderson2} for some background on this subject.)

Two preliminary observations should be made.
First, since the Lorentzian norm of a non-zero tensor may vanish
it is clear that only limited information would be gained from an
assumption on the Lorentzian norm of the curvature tensor.
This justifies that we endow the Lorentzian manifold with a ``reference'' Riemannian metric
(denoted by $g_T$ below); this metric is defined 
at a point $p$ once we prescribe a future-oriented time-like unit vector $T$. 
We refer to the pair $(p,T)$ as an observer located at the point $p$. 
This reference vector is necessary in order to define appropriate notions of conjugate and injectivity radii. 
(See Section~\ref{LO-0} below for details.)

Secondly, we rely here on the elementary but essential observation that, in the flat 
Riemannian and Lorentzian spaces, geodesics (are straight lines and therefore) coincide.
Under our assumptions, we will see that geodesics associated with the given 
Lorentzian metric are comparable to geodesics associated with the reference Riemannian
metric. On the other hand, the curvature bound assumed on the Lorentzian metric 
implies, in general, no information on the curvature of the reference metric. As we show below, 
one of the main issues is to guarantee the regularity of a foliation of the manifold by spacelike hypersurfaces.


\subsection*{Earlier work}

Let us briefly review some classical results from Riemannian geometry. Let $(M,g)$ be a
differentiable $n$-manifold (possibly with boundary) endowed with a 
Riemannian metric $g$. (Throughout the present paper, the manifolds and metrics are always assumed
to be smooth.)
Denote by $\Bcal(p,r)$ the corresponding geodesic ball centered at $p \in M$ and with radius $r>0$.
Suppose that at some point $p\in M$ the unit ball $\Bcal(p,1)$ is compactly contained in $M$
and that the Riemann curvature bound and the lower volume bound, 
\be
\label{Linfinity} \|\Riem_g\|_{\Lbf^\infty(\Bcal(p,1))}\leq K,
\qquad
\vol_g(\Bcal(p,1)) \geq v_0, 
\ee
hold for some constants $K, v_0 >0$.
(We use the standard notation $\Lbf^m$, $1 \leq m \leq \infty$, for the spaces of Lebesgue measurable functions.)
Then, according to Cheeger, Gromov, and Taylor \cite{CGT}
the injectivity radius $\Inj_g(M,p)$ at the point $p$ is bounded below by a positive constant $i_1=i_1(K,v_0,n)$,
\be
\label{inj1}
\Inj_g(M,p) \geq i_1.
\ee
It should be noticed that this is a local statement; for earlier work on the injectivity radius 
see \cite{Cheeger1,CLY,HeintzeKarcher}.

Moreover, Jost and Karcher \cite{JostKarcher} relied on the regularity theory for elliptic operators
and established the existence of coordinates in which the metric has optimal regularity
and are defined in a ball with radius $i_2=i_2(K, v_0,n)$. 
Precisely,  given $\eps >0$ and $0<\gamma<1$ there exist
a positive constant $C(\eps,\gamma)$ (depending also upon $(K,v_0,n$)
and a system of harmonic coordinates defined in the geodesic ball $\Bcal(p, i_2)$
in which the metric $g$ is close to the Euclidian metric $g_E$ in these coordinates 
and has optimal regularity, in the following sense:
\be
\label{epsilon}
\aligned
& e^{-\eps} \, g_E \leq g \leq  e^\eps \, g_E, 
\\
& r^{1+\gamma}\| \del g \|_{\Cbf^{\gamma}(\Bcal(p,r))} \leq C(\eps,\gamma),
\qquad r \in (0,i_2].
\endaligned
\ee
Here, $\Cbf^0$ and $\Cbf^{\gamma}$ are
the spaces of continuous and H\"older continuous functions, respectively.
Harmonic coordinates are optimal \cite{DeTurckKazdan} in the sense that
if the metric is of class $\Cbf^{k,\gamma}$ in certain coordinates then it has at least the same
regularity in harmonic coordinates.

The above results were later generalized by Anderson \cite{Anderson0} 
and Petersen \cite{Petersen} who replaced the $\Lbf^\infty$ curvature bound
by an $\Lbf^m$ curvature bound with $m>n/2$. 
For instance, one can take $m=2$ in dimension $n=3$ in the application to general relativity 
(since time-slices of Lorentzian $4$-manifolds are Riemannian $3$-manifolds).

It is only more recently that the same questions were tackled for
Lorentzian $(n+1)$-manifolds $(M,g)$.
Anderson \cite{Anderson1,Anderson2} studied the long-time evolution of solutions to the Einstein field 
equations and formulate several conjecture. 
In particular, assuming the Riemann curvature bound in some domain $\Omega$
\be
\label{Linfinity2}
\|\Riem_g\|_{\Lbf^\infty(\Omega)} \leq K
\ee
and other regularity conditions, he investigated the existence of coordinates
that are harmonic in each spacelike slice of a time foliation of $M$.
This work by Anderson motivated our investigation in the present paper.

On the other hand, motivated by applications to general relativity and
nonlinear wave equations using harmonic analysis tools,
Klainerman and Rodnianski \cite{KR4} considered asymptotically flat spacetimes
endowed with a time foliation and satisfying the $\Lbf^2$ curvature bound
\be
\label{Lsquare2}
\|\Riem_g\|_{\Lbf^2(\Sigma)} \leq K
\ee
for every spacelike hypersurface $\Sigma$. 
To control the injectivity radius of past null cones, they relied on their earlier work
\cite{KR1,KR2} on the conjugate radius of null cones in terms of Bell-Robinson's
energy and energy flux, and derived in \cite{KR4} a new estimate for the null cut locus radius.
We refer to these papers for further details and references on the Einstein equations. Section~\ref{nullsection} of
the present paper is a prolongation of the work \cite{KR4}.

\

\subsection*{Outline of this paper}

The present paper establishes four estimates on the radius of injectivity of Lorentzian
manifolds, which hold either in a neighborhood of a point or on the past null cone at a point.
Our assumptions are formulated within a geodesic ball (or within a null cone)
and possibly apply in a ball with arbitrary size as long as our curvature and volume assumptions hold.
All assumptions and statements are local and geometric.

An outline of the paper is as follows.
In Section~\ref{LO-0}, we begin with basic material from Lorentzian geometry and we introduce
the notions of reference metric and exponential map for Lorentzian manifolds.
In Section~\ref{WS-0}, we state our first estimate (Theorem~\ref{inject} below) for a class of manifolds
that have bounded curvature and admit a time foliation by slices with bounded extrinsic curvature.
In Section~\ref{proofinject}, we provide a proof of this first estimate and we introduce a technique
that will be used (in variants) throughout this paper;
we combine two main ingredients : sharp estimates for Jacobi fields along geodesics, and
an homotopy argument based on contracting a possible loop to two linear segments.
In Section~\ref{convexsection}, our second main result (Theorem~\ref{convex}) shows,
under the same assumptions, 
the existence of convex functions (distance functions) and convex neighborhoods; this result leads us
to a lower bound of the convexity radius.

In Section~\ref{nullsection}, our third estimate (Theorem~\ref{nulltheorem}) covers the case of null cones
under the assumption that the manifold has $L^2$ bounded curvature
on every spacelike slice; this provides a generalization and
an alternative proof to the result by Klainerman and Rodnianski in \cite{KR4}.

Next, in Section~\ref{nopre}, we establish our principal and fourth result (stated in Theorem~\ref{nofoliation}
below) which provides an injectivity radius bound under the mild assumption that the exponential map $\expb_p$
is defined in some ball and the curvature $\Riem$ is bounded.
Most importantly, this is a general result that does not require a time foliation of the manifold
but solely a single reference (future-oriented time-like unit) vector $T$ at the base point $p$.
This is very natural in the context of general relativity and $(p,T)$ is interpreted as an observer
at the point $p$.

Given an observer $(p,T)$, we can define the ball $B_T(0,r) \subset T_pM$
with radius $r$, determined by the reference Riemannian inner product at $p$,
and we can also define the geodesic ball $\Bcal_T(p,r):=\expb_p(B_T(0,r))$.
In turn, the radius of injectivity $\Inj_g(M,p,T)$ is defined as the largest radius
$r$ such that the exponential map is a diffeomorphism from $B_T(0,r)$ onto $\Bcal_T(p,r)$.
Let us then consider an arbitrary geodesic $\gamma=\gamma(s)$ 
initiating at $p$ and let us $g$-parallel transport the vector $T$ along this geodesic, 
defining therefore a vector field $T_\gamma$ along this geodesic, only. At every point of $\gamma$ 
we introduce the reference metric $g_{T_\gamma}$ and compute the curvature norm $|\Riem_g|_{T_\gamma}$. 
This allows us to express the curvature bound. 
For the convenience of the reader we state here our main result and refer to Section~\ref{nopre}
for further details.

\begin{theorem}[Injectivity radius of Lorentzian manifolds]
\label{nofoliation}
Let $M$ be a time-orientable Lorentzian, differentiable $(n+1)$-manifold.
Consider an observer $(p,T)$ consisting of a point $p\in M$ and a reference 
(future-oriented time-like unit) vector $T\in T_pM$. 
Assume that the exponential map $\expb_p$ is defined in the ball $B_T(0,r_0) \subset T_pM$
and the Riemann curvature satisfies
\be
\label{hypot7}
\sup_\gamma | \Riem_g|_{T_\gamma} \leq {1 \over r_0^2}, 
\ee
where the supremum is over the domain of definition of $\gamma$ and 
over every $g$-geodesic $\gamma$ initiating from a vector in the Riemannian ball $B_T(0,r_0) \subset T_pM$.
Then, there exists a constant $c(n)$  depending only on the dimension $n$ such that
\be
\label{inequa1}
{\Inj_g(M,p,T) \over r_0} \geq c(n) \, \frac{\vol_g\big(\Bcal_T(p,c(n)r_0)\big)}{r_0^{n+1}}.
\ee
\end{theorem}

This result should be compared with the injectivity radius estimate established by Cheeger, Gromov, and 
Taylor \cite{CGT} in Riemannian geometry. Observe that the curvature assumption \eqref{hypot7} can 
always be satisfied by suitably rescaling the metric tensor. 
It would be interesting to replace the volume term in the right-hand side 
of \eqref{inequa1} by $\vol_g\big(\Bcal(p,r_0)\big)$.

Finally, in the last two sections of this paper, 
we establish a volume comparison theorem for future cones
and generalize our main theorem to the volume of a future cone (Section~\ref{conclude1}), and 
we briefly discuss the regularity of the Lorentzian metric in harmonic-like coordinates, 
and present a generalization to pseudo-riemannian manifolds (Section~\ref{conclude2}).


\section{Preliminaries on Lorentzian geometry}
\label{LO-0}

\subsection*{Basic definitions}
It is useful to discuss first some basic definitions from Lorentzian geometry,
for which we can refer to the textbook by Penrose \cite{Penrose}.
Throughout this paper, $(M,g)$ is a connected 
and differentiable $(n+1)$-manifold, endowed with a Lorentzian metric tensor
$g$ with signature $(-, +, \ldots, +)$.
To emphasize the role of the metric $g$ or the point $p$ we use any of the following notations
$$
g_p(X,Y) = \langle X, Y \rangle_{g_p} = \langle X, Y \rangle_g = \langle X, Y \rangle_p
$$
for the inner product of two vectors $X,Y$ at a point $p \in M$;
we sometimes also write $|X|_{g_p}^2$ instead of $g_p(X,X)$.
Recall that the tangent vectors $X \in T_p M$ are called time-like, null, or spacelike
depending whether the norm $g_p(X,X)$ is negative, zero, or positive, respectively.
Vectors that are time-like or null are called causal.

The time-like vectors form a cone with two connected components.
The manifold $(M,g)$ is said to be time-orientable if we can select in a continuous way
a half-cone of time-like vectors at every point $p$. The choice of a specific orientation allows
us to decompose the cone of time-like vectors into future-oriented and past-oriented ones.
The set of all future-oriented, time-like vectors at $p$ and the corresponding bundle on $M$
are denoted by $T^+_pM$ and $T^+M$, respectively. We also introduce the bundle $T_1^+M$ 
consisting of elements of $T^+M$ with unit length. 

By definition, a trip is a continuous curve $\gamma : (a,b) \to M$ made of finitely many
future-oriented, time-like geodesics. We write $p << q$ if there exists a trip from $p$ to $q$.
A causal trip is defined similarly except that the geodesics may be causal instead of time-like,
and we write $p < q$ if there exists a causal trip from $p$ to $q$. 

The set $\Ical^+(p) := \big\{ q \in M \, / \, p << q \big\}$ is called the chronological future of the point $p$,
and
$\Ical^-(p) := \big\{ q \in M \, / \, q << p \big\}$ is called the chronological past. The causal future and past
are defined similarly by replacing $<<$ by $<$. The future or past sets of a set $S \subset M$ are defined by
$$
\Ical^\pm(S) := \bigcup_{p \in S} \Ical^\pm(p),  \qquad  \Jcal^\pm(S) := \bigcup_{p \in S} \Jcal^\pm(p),
$$
and one easily checks that $\Ical^\pm(S)$ are open, but that $J^\pm(S)$ need not be closed in general.

A future set $F \subset M$ by definition has the form $F = \Ical^+(S)$ for some set $S \subset M$. Similarly,
a past set satisfies $F = \Ical^-(S)$ for some $S$. A set is called achronal if no two points are connected
by a time-like trip. Observe that a set can be spacelike at every point without being achronal
and that an achronal set can be null at some (or even at every) point.
A set $B \subset M$ is called an achronal boundary if it is the boundary of a future set, that is
$B = \del \Ical^+(S) = \overline{\Ical^+(S)} \setminus \Ical^+(S)$ for some $S \subset M$. One can check that
given a non-empty achronal boundary the manifold can be partitioned as $M = P \cup B \cup F$, where
$B$ is the boundary of both $F$ and $P$ and, moreover, any trip from $p \in P$ to $q \in F$ meets $B$
at a unique point. Observe also that any achronal boundary is a Lipschitz continuous $n$-manifold.
For instance, in Section~\ref{nullsection} below, we will be interested in the geometry of past null cones,
that is the sets $\del \Jcal^-(p)$ for $p \in M$.

Given an arbitrary achronal and closed set $S \subset M$,
we define the (future or past) domains of dependence of $S$ in $M$ by
$$
\aligned
& \Dcal^\pm(S) := \big\{ p \in M \, / \, \text{every future (resp. past) endless trip containing $p$ meets
$S$} \big\},
\\
& \Dcal(S) := \Dcal^-(S) \cup \Dcal^+(S).
\endaligned
$$
Observe that domains of dependence are closed sets. Next, define the (future or past) Cauchy horizons
$$
\aligned
& \Hcal^\pm(S) := \big\{ p \in \Dcal^\pm(S) \, / \, \Ical^\pm(p) \cap \Dcal^\pm(S) = \emptyset \big\}
                     = \Dcal^\pm(S) \setminus \Ical^\mp(\Dcal^\pm(S)),
\\
& \Hcal(S) := \Hcal^-(S) \cup \Hcal^+(S).
\endaligned
$$
For instance, the future Cauchy horizon is the future boundary of the future domain of dependence of $S$.
One can check that the Cauchy horizons are closed and achronal sets, with $\del \Dcal^+(S) = \Hcal^+(S) \cup S$
and $\del \Dcal(S) = \Hcal(S)$.

Finally, a (future) Cauchy hypersurface for $M$ is defined as an achronal (but not necessarily closed) set $S$
satisfying $\Dcal^+(S)=M$. For instance, it is sufficient for $\overline{S}$ to be smooth, achronal, spacelike
and such that every endless null geodesic meet $M$.


\subsection*{Reference metric}

As explained in the introduction one should not use the Lorentzian metric to compute the norm of a tensor
since the Lorentzian norm may vanish even when the tensor does not. This motivates the introduction
of a ``reference'' Riemannian metric associated with a time-like vector field, as follows.

Let $T$ be a future-oriented, time-like, unit vector field, satisfying therefore
${g_p(T,T) = -1}$ at every point $p$. We refer to $T$ as the {\sl reference vector field}
prescribed on $M$. Introduce a moving frame $E_\alpha$ ($\alpha =0,1, \ldots, n$)
defined in $M$,
that is, $E_\alpha$ is an orthonormal basis of the tangent space at every point
and consists of the vector $e_0=T$ supplemented with $n$ spacelike unit vectors $e_j$ $(j=1, \ldots, n$).
Denoting by $E^\alpha$ the corresponding dual frame, the Lorentzian metric tensor takes the form
$$
g = \eta_{\alpha \beta} \, E^\alpha \otimes E^\beta,
$$
where $\eta_{\alpha\beta}$ is the Minkowski ``metric''. This decomposition suggests to consider the Riemannian version
obtained by switching the minus sign in $\eta_{00}=-1$ into a plus sign, that is
$$
g_T := \delta_{\alpha \beta} \, E^\alpha \otimes E^\beta,
$$
where $\delta_{\alpha \beta}$ is the Euclidian ``metric".
Clearly, $g_T$ is a positive definite metric; it is referred to as the {\sl reference Riemannian metric}
associated with the frame $E_\alpha$.

For every $p \in M$, since $T_p$ is time-like,
the restriction of the metric $g_p$ to the orthogonal complement
$\{T_p \}^\perp \subset T_PM$ is positive definite, and
the reference metric can be computed as follows:
if $V = a \, T_p + V'$ and $W = b \, T_p + W'$ with $V',W' \in \{T_p \}^\perp$, then
$$
g_{T,p}(V,W) = a \, b + g_p(V',W').
$$
In the following, we use the notation
$$
g_{T,p}(V,W) = \langle V,W \rangle_{T,p},
\qquad
g_{T,p}(V,V) = |V|_{T,p}^2
$$
for vectors; the norm of tensors is defined and denoted similarly.

In contrast with the Lorentzian norm, the Riemannian norm  $|A|_{T,p}$ of a tensor $A$ at a point $p \in M$
vanishes
if and only if the tensor is zero at $p$. Moreover, as long as $T$ remains in a compact subset
of the bundle of half-cone $T^+M$, the norms associated with different reference vectors are equivalent.

The reference Riemannian metric also allows one to define the functional norms for
Lebesgue and Sobolev spaces of tensors defined on $M$ (as well as on submanifolds of $M$),
allowing us for instance to view $\Lbf^2(M,g_T)$ as a Banach space.
In particular, we will use later the $\Lbf^2$ norm of a tensor field $T$ on $M$
restricted to an hypersurface $\Sigma$:
$$
\| \nabla h \|_{\Lbf^2(\Sigma, g_T)}  := \int_\Sigma |\nabla h|_T^2 \, dV_{\Sigma,g_T},
$$
where $dV_{\Sigma,g_T}$ is the volume form induced on $\Sigma$ by the reference Riemannian metric.
The functional norm above depends upon the choice of the vector field $T$, but another choice of $T$
would give rise
to an equivalent norm (provided $T$ remains in a fixed compact subset).
Observe in passing that the volume forms associated with the metrics $g$ and $g_T$ coincide,
so that the spacetime integrals of {\sl functions} in $(M,g)$ or $(M,g_T)$ coincide;
for instance, the volume $\vol_g(A)$ and $\vol_{g_T}(A)$ of a set $A \subset M$ coincide.

Furthermore, we observe that in order to define the reference metric $g_T$ at a given point $p$,
it suffices to prescribe a future-oriented time-like unit {\sl vector} $T$ at that point $p$ only;
it is not necessary to prescribe a vector field.
In the situation where the reference metric need only be defined at the base point $p$,
we refer to $T$ as the {\sl reference vector} (rather than vector field)
and we refer to $(p,T) \in T_1^+M$ as the {\sl observer at the point $p$.}
This will be the standpoint adopted for our main result in Section~\ref{nopre} below.

\subsection*{Exponential map}

On a complete Riemannian manifold the exponential map $\expb_p : T_pM \to M$ at some point $p \in M$
is defined on the whole tangent space $T_pM$ and is smooth.
For sufficiently small radius $r$ the restriction of $\expb_p$ to the ball $B_{g_p}(0,r) \subset T_pM$
(determined by the metric $g_p$ at the point $p$)
is a diffeomorphism on its image. The radius of injectivity at the point $p$ is defined as the largest value of $r$
such that the restriction $\expb_p|_{B_{g_p}(0,r)}$ is a global diffeomorphism.

In the Lorentzian case, the exponential map is defined similarly but
some care is needed in defining the notion of radius of injectivity. First of all, if the manifold is not
geodesically complete (which is a rather generic situation, as illustrated by
Penrose and Hawking's incompleteness theorems \cite{HawkingEllis}),
the map $\expb_p$ need not be defined on the whole tangent space $T_pM$
but only on a neighborhood of the origin in $T_pM$. More importantly,
the Lorentzian norm of a non-zero vector may well vanish;
consequently, the radius of injectivity should not be defined directly from the Lorentzian metric $g$.
The definition below depends on the prescribed Riemannian metric $g_{T,p}$ at the point $p$.

\begin{definition}
\label{radius}
The {\sl conjugate radius} $\Conj_g(M,p,T)$ of an observer $(p,T) \in T_1^+M$ 
is the largest radius $r$ such that the exponential map $\expb_p$
is a local diffeomorphism from the Riemannian ball $B_T(0,r) = B_{g_{T,p}}(0,r) \subset T_pM$
to a neighborhood of $p$ in the manifold $M$.
Similarly, the {\rm injectivity radius} $\Inj_g(M,p,T)$ of an observer $(p,T) \in T_1^+M$ 
is the largest radius $r$ ssuch that
the exponential map is a global diffeomorphism at every point of the ball $B_T(0,r)$.
\end{definition}

When a vector field $T$ is prescribed on the manifold (rather than a vector at the point $p$),
we use the notation $\Inj_g(M,p,T)$ instead of $\Inj_g(M,p,T_p)$.
The radii $\Conj_g(M,p,T)$ and $\Inj_g(M,p,T)$ are essentially independent of
the choice of the reference vector, as long as it remains in a fixed compact subset of $T_p^+M$.

We also need the notion of injectivity radius for null cones. Given a point $p \in M$ and
a reference vector $T \in T_p M$, we consider the {\sl past null cone} at $p$,
$$
N_p^- := \big\{ X \in T_pM \, \big/ \, g_p(X,X) = 0, \, g_p(T, X) \geq 0 \big\},
$$
which is defined a subset of the tangent space at $p$.
Denote by
$$
B_T^N(0,r) = B_{g_{T,p}}^N(0,r) := B_{g_{T,p}}(0,r)\cap N_p^-
$$
the intersection of the Riemannian $g_{T,p}$-ball with radius $r$ and the past null cone, and by
$$
\Ncal^-(p) := \del \Ical^-(p)
$$
the past null cone at $p$.

Consider now the restriction of $\expb_p$ to the past null cone, denoted by
$$
\expb^N_p : B_T^N(0,r) \subset N_p^- \to \Ncal^-(p) \subset M,
$$
which we refer to as the {\sl null exponential map.}

\begin{definition}
\label{nullradius} 
The {\rm (past) null conjugate radius} $\NullConj_g(M,p,T)$
of an observer $(p,T) \in T_1^+M$ 
is the largest radius $r$ such that the null exponential map $\expb_p^N$ is a local diffeomorphism
from the punctured Riemannian ball $B_T^N(0,r) \backslash\{0\}\subset T_pM$
to a neighborhood of $p$ in the past null cone. The {\rm null injectivity radius} $\NullInj_g(M,p,T)$
of an observer $(p,T) \in T_1^+M$ 
is defined similarly by requiring the map $\expb_p^N$ to be a global diffeomorphism.
\end{definition}


\section{Lorentzian manifold endowed with a reference vector field}
\label{WS-0}

\subsection*{A first injectivity radius estimate}

From now on, we fix a reference vector field $T$ which allows us to define the Riemannian metric $g_T$
and compute the norms of tensors.
We begin with a set of assumptions encompassing a large class of Lorentzian manifolds
with $\Lbf^\infty$ bounded curvature and we state our first injectivity estimate, in Theorem~\ref{inject} below.
The forthcoming sections will be devoted to further generalizations and variants of this result.

We fix a point $p\in M$ and assume that a domain $\Omega \subset M$ containing $p$
is foliated by spacelike hypersurfaces $\Sigma_t$ with future-oriented time-like unit normal~$T$,
\be
\label{foliation1}
\Omega = \bigcup_{t \in [-1,1]} \Sigma_t.
\ee
The positive coefficient $n$ is defined by the relation ${\del \over \del t} = n \, T$, or
$$
n^2 := - g\big({\del \over \del t}, {\del \over \del t} \Big).
$$
In the context of general relativity, $n$ is the proper time of an observer moving orthogonally
to the hypersurfaces, and is called the {\sl lapse function}. The geometry of the foliation
is determined by this function $n$ together with the Lie derivative $\Lie_Tg$. The latter
is nothing but the second fundamental form, or extrinsic curvature, of the slices $\Sigma_t$
embedded in the manifold $M$.

We always assume that the geodesic ball $\Bcal_{\Sigma_0}(p,1) \subset \Sigma_0$
(determined by the induced metric $g|_{\Sigma_0}$) is compactly contained in $\Sigma_0$.
We introduce the following assumptions:
\begin{equation*}
\tag{A1}
e^{-K_0} \leq n \leq e^{K_0} \quad \text{ in } \Omega, 
\end{equation*}
\begin{equation*}
\tag{A2}
|\Lie_Tg |_T \leq K_1 \quad \text{ in } \Omega,
\end{equation*}
\begin{equation*}
\tag{A3} |\Riem_g |_T\leq K_2 \quad  \text{ in } \Omega,
\end{equation*}
\begin{equation*}
\tag{A4} \vol_{g|_{\Sigma_0}}(\Bcal_{\Sigma_0}(p,1))\geq v_0,
\end{equation*}
where $K_0, K_1, K_2$ and $v_0$ are positive constants.
Observe that Assumption~$(A4)$ is a condition on the initial slice only; together with the other assumptions
it actually implies a lower volume bound for every slice of the foliation.

We will prove:

\begin{theorem}[Injectivity radius of foliated manifolds]
\label{inject}
Let $M$ be a 
differentiable manifold endowed with a Lorentzian metric $g$
satisfying the regularity assumptions $(A1)$--$(A4)$ at some point $p$ and for some foliation \eqref{foliation1}.
Then, there exists a positive constant $i_0$ depending only upon the foliation bounds $K_0,K_1$, the curvature bound $K_2$, the volume bound $v_0$, and the dimension
of the manifold such that the injectivity radius at $p$ satisfies
$$
\Inj_g(M,p,T) \geq i_0.
$$
\end{theorem}

The following section is devoted to the proof of this theorem. Observe that
the conditions $(A1)$--$(A4)$ are {\sl local} about one point of the manifold
and are stated in purely geometric terms, requiring no particular choice of coordinates.
Of course, the conclusion of Theorem~\ref{inject} hold globally in $M$
if the assumptions $(A1)$--$(A4)$ hold also globally at every point of the manifold.
Our assumptions do depend on the choice of the time-like vector field $T$, but the dependence
of the constants arising in $(A1)$--$(A4)$ should not be essential; however, it is conceivable that,
when applying this theorem in a specific situation a quantitatively sharper estimate would be obtained
with a choice of an ``almost Killing'' field, that is a field $T$ corresponding to a ``small''
Lie derivative $\Lie_T g$. Later in Section~\ref{nopre}, 
a more general approach is presented in which the vector field $T$ is constructed 
from a single vector prescribed at the point $p$. 


\subsection*{Basic estimates on the reference metric}

To establish Theorem~\ref{inject} it is convenient to introduce coordinates on $\Omega$, 
chosen as follows.
Fix arbitrarily some coordinates $(x^i)$  on the initial slice $\Sigma_0$.
Then, transport these coordinates to the whole of $\Omega$ along the integral curves of the vector field $T$.
This construction generates coordinates $(x^\alpha)$ on $\Omega$ such that
$x^0=t$ and the vector $\del/\del t$ is orthogonal to $\del/\del x^j$, so that the Lorentzian metric takes the form
\be
\label{foliation2} g =- n^2 \, dt^2 + g_{ij} dx^i dx^j,
\ee
where $n$ is the lapse function and $g_{ij}$ is the Riemannian metric induced on the slices $\Sigma_t$.
The reference Riemannian metric in the domain $\Omega$ then takes the form
\be
\label{refer}
g_T = n^2 dt^2 + g_{ij} dx^i dx^j,
\ee
and the Riemannian norm of a vector $X$ has the explicit form: $g_T(X,X) := n^2 \, X^0 X^0 + X^j X_j$.

We want to control the discrepancy
between the reference Riemannian metric $g_T$ and the original Lorentzian metric $g$,
as measured in the connections $\nabla$ and $\nabla_{g_T}$ and the 
curvature tensors $\Riem$ and $\Riem_{g_T}$. Clearly,
these estimates should involve the constants arising in $(A1)$--$(A4)$.
Consider the general class of metrics
\be
\label{general}
\gt: = f \, dt^2 + g_{ij} \, dx^idx^j,
\ee
which allows us to recover both the Lorentzian ($f=-n^2$) and the Riemannian ($f=n^2$) metrics.

In view of the expressions of the Christoffel symbols and the Riemann curvature
\be
\nonumber
\begin{aligned}
& \Gammat^\gamma_{\alpha\beta}
   : = {1 \over 2} \gt^{\gamma\delta} \, \big( \frac{\del \gt_{\delta\beta}}{\del x^\alpha}+\frac{\del
        \gt_{\delta\alpha}}{\del x^{\beta}}-\frac{\del \gt_{\alpha\beta}}{\del x^{\delta}}\big),
\\
& \Rt^{\zeta}_{\alpha\beta\delta}
      := \frac{\del \Gammat^{\zeta}_{\beta\delta}}{\del  x^\alpha}-\frac{\del \Gammat^{\zeta}_{\alpha\delta}}
           {\del x^{\beta}}+\Gammat^{\zeta}_{\alpha\eta}\Gammat^{\eta}_{\beta\delta}
         -\Gammat^{\zeta}_{\beta\eta}\Gammat^{\eta}_{\alpha\delta},
\\
& \Rt_{\alpha\beta\gamma\delta} := \gt_{\gamma\zeta}\Rt^{\zeta}_{\alpha\beta\delta},
\qquad
\Rt_{\alpha\beta} := \Rt_{\alpha\gamma\beta\delta}\gt^{\gamma\delta},
\end{aligned}
\ee
we compute explicitly the Christoffel symbols associated with the metric $\gt$,
\be
\label{chris}
\aligned
& \Gammat^0_{00} =\frac{1}{2 f}\frac{\del f}{\del t},
    \qquad
   \Gammat^0_{0i} =\frac{1}{2 f}\frac{\del f}{\del x^i},
   \qquad  \Gammat^0_{ij} =-\frac{1}{2 f}\frac{\del g_{ij}}{\del t},
\\
&  \Gammat^{k}_{00} =-\frac{1}{2 }g^{kl}\frac{\del f }{\del x^{l}},
     \qquad \Gammat^{k}_{i0} =\frac{1}{2 }g^{kl}\frac{\del g_{li} }{\del t},
     \qquad \Gammat^{k}_{ij} ={\Gamma}^{k}_{ij},
\endaligned
\ee
as well as the non-trivial curvature terms
$$
\Rt_{ijkl}=R_{ijkl}-\frac{1}{4f} \, \big( \frac{\del g_{ik}}{\del t}\frac{\del g_{jl}}{\del
t}-\frac{\del g_{il}}{\del t}\frac{\del g_{jk}}{\del t}\big),
$$
\begin{equation*}
\begin{split}
\Rt^{p}_{0jl} 
& =\frac{\del {\Gamma}^{p}_{jl}}{\del t}
    - {\del \over \del x^j}
     \, \big( {1 \over 2} g^{pq}\frac{\del g_{ql}}{\del t}\big)
     +{1 \over 2} g^{pq}\frac{\del}{\del t}
g_{qk}{\Gamma}^{k}_{jl} - {\Gamma}^{p}_{jk}({1 \over 2} g^{kq}\frac{\del g_{lq}}{\del t})
\\
& \quad + \big(\frac{1}{4f}g^{pq}\frac{\del f}{\del x^q}\big) \, \frac{\del g_{jl}}{\del t}
-{1 \over 2} g^{pq}\frac{\del g_{qj}}{\del t}\frac{1}{2f}\frac{\del f}{\del x^{l}},
\end{split}
\end{equation*}
$$
\Rt_{0jil} = {1 \over 2}
     \, \big( \nabla_{l} \big( \frac{\del}{\del t}g_{ij}\big) - \nabla_i \big( \frac{\del}{\del t}g_{lj}\big)\big)
    +\frac{1}{4f} \, \big( \frac{\del f}{\del x^i}\frac{\del g_{jl}}{\del t}
    -\frac{\del f}{\del x^{l}}\frac{\del g_{ij}}{\del t} \big),
$$ 
\begin{equation*}
\begin{split}
\Rt^{p}_{i00}
& = \frac{\del }{\del x^i} \big( -{1 \over 2} g^{pq}\frac{\del f}{\del x^q} \big)
   -\frac{\del } {\del t} \big( {1 \over 2} g^{pq}\frac{\del g_{qi}}{\del t} \big)
   + {\Gamma}^{p}_{il} \big( -{1 \over 2} g^{lq}\frac{\del f}{\del x^q} \big)
\\
& \ \ \ - \big({1 \over 2} g^{pq}\frac{\del g_{ql}}{\del t}\big) \,
     \big({1 \over 2} g^{lr}\frac{\del g_{ri}}{\del t}\big) + \frac{1}{2f} \frac{\del f}{\del t}
       \big( {1 \over 2} g^{pq}\frac{\del g_{qi}}{\del t} \big)
      + {1 \over 2} g^{pq}\frac{\del f}{\del x^q}\frac{1}{2f}\frac{\del f}{\del x^i},
\end{split}
\end{equation*}
and
\begin{equation*}
\begin{split}
\Rt_{i0j0}
& = - {1 \over 2} \big(\nabla_i\nabla_{j}f+\frac{\del^2 g_{ij}}{\del t^2}\big)
    + \frac{1}{4}g^{pq}\frac{\del g_{ip}}{\del t}\frac{\del g_{jq}}{\del t}
    +\frac{1}{4f}\frac{\del f}{\del t}\frac{\del g_{ij}}{\del t}
    +\frac{1}{4f}\frac{\del f}{\del x^i}\frac{\del f}{\del x^j}.
\end{split}
\end{equation*}

By applying the formulas above to both metrics $g, g_T$ we estimate the Christoffel symbols, as follows.
Recall that the difference $\Gamma^\gamma_{\alpha\beta} - \Gamma^\gamma_{g_T,\alpha\beta}$
can be regarded as a tensor field on $M$, so
that the following (Riemannian) norm squared
is a scalar field on the manifold $M$:
$$
| \nabla_{g_T} - \nabla|_T^2
:= | \Gamma_{g_T} - \Gamma|_T^2
=   (\Gamma^\alpha_{g_T,\beta\gamma}      - \Gamma^\alpha_{\beta\gamma}) \,
     (\Gamma^{\alpha'}_{g_T,\beta'\gamma'} - \Gamma^{\alpha'}_{\beta'\gamma'}) \,
     g_{T,\alpha\alpha'}\, g_T^{\beta\beta'} \, g_T^{\gamma\gamma'}.
$$
We need also the expression of the Lie derivative of $g$ along the vector field $T$.
By a direct computation from \eqref{foliation2} we obtain
\be
\label{Lie}
(\Lie_Tg)_{00} =0,
\qquad
(\Lie_Tg)_{0i} = \frac{1}{n}\frac{\del n}{\del x^i},
\qquad
(\Lie_Tg)_{ij}= \frac{1}{n}\frac{\del g_{ij}}{\del t}.
\ee

\begin{lemma}[Levi-Cevita connection of the reference metric]
\label{Gamma}
Suppose that $g$ satisfies Assumptions $(A1)$-$(A2)$.
Then, the covariant derivative of the Lorentzian and Riemannian metrics
are comparable, precisely  
$$
| \nabla_{g_T} - \nabla|_T = n^2 \, |\Lie_T g|_T^2 
\leq e^{2K_0} K_1^2 =:K_3. 
$$
\end{lemma}

\begin{proof} In view of \eqref{chris} the difference $\Gamma_{g_T} - \Gamma$ depends
essentially upon the terms $\frac{\del n}{\del x^i}$ and $\frac{\del g_{ij}}{\del t}$
which precisely appear in the expression of the Lie derivative \eqref{Lie}.
We omit the details.
\end{proof}

It is important to observe that the difference between the curvature tensors 
can not be similarly estimated, and that this is one of the main difficulties
to deal with in the present work. 

For future reference we provide here the expressions of
certain curvature coefficients of $g$ and $g_T$ in terms of (first-order derivatives of)
the lapse function $n$ and the induced metric $g_{jk}$:
\begin{equation*}
\begin{split}
R_{ijkl}
& = R^{\Sigma}_{ijkl}+\frac{1}{4n^2} \, \big( \frac{\del g_{ik}}{\del t}\frac{\del g_{jl}}{\del t}
   -\frac{\del g_{il}}{\del t}\frac{\del g_{jk}}{\del t}\big),
\\
R_{0jil}
& = {1 \over 2}  \big( \nabla_{l}(\frac{\del}{\del t}g_{ij} \big)
     -\nabla_i \big(\frac{\del}{\del t}g_{lj}) \big)
     + \frac{1}{4n^2} \big( \frac{\del n^2}{\del x^i}\frac{\del g_{jl}}{\del t}
     -\frac{\del n^2}{\del x^{l}}\frac{\del g_{ij}}{\del t} \big),
\\
R_{i0j0}
& =   {1 \over 2} \big (\nabla_i\nabla_{j}(n^2)-\frac{\del^2 g_{ij}}{\del t^2}\big)
     + \frac{1}{4}g^{pq}\frac{\del g_{ip}}{\del t}\frac{\del g_{jq}}{\del t}
     + \frac{1}{4n^2}\frac{\del n^2}{\del t}\frac{\del g_{ij}}{\del t}
     - \frac{1}{4n^2}\frac{\del n^2}{\del x^i}\frac{\del n^2}{\del x^j},
\end{split}
\end{equation*}
and
\begin{equation*}
\begin{split}
R_{T,ijkl}
& = R_{ijkl}^{\Sigma}-\frac{1}{4n^2} \big( \frac{\del g_{ik}}{\del t}\frac{\del g_{jl}}{\del t}
    - \frac{\del g_{il}}{\del t}\frac{\del g_{jk}}{\del t} \big),
\\
R_{T,0jil}
& =  {1 \over 2}  \big( \nabla_{l} \big(\frac{\del}{\del t}g_{ij}\big)
     -\nabla_i \big( \frac{\del}{\del t}g_{lj}) \big)
     + \frac{1}{4n^2} \big( \frac{\del n^2}{\del x^i}\frac{\del g_{jl}}{\del t}
     - \frac{\del n^2}{\del x^{l}}\frac{\del g_{ij}}{\del t} \big),
\\
R_{T,i0j0}
& = {1 \over 2} \big( \nabla_i\nabla_{j}(-n^2)-\frac{\del^2 g_{ij}}{\del t^2}\big)
    + \frac{1}{4}g^{pq}\frac{\del g_{ip}}{\del t}\frac{\del g_{jq}}{\del t}
    +\frac{1}{4n^2}\frac{\del n^2}{\del t}\frac{\del g_{ij}}{\del t}
    + \frac{1}{4n^2}\frac{\del n^2}{\del x^i}\frac{\del n^2}{\del x^j},
\end{split}
\end{equation*}
where $R_{ijkl}^{\Sigma}$ denotes the induced curvature tensor on the time slices $\Sigma = \Sigma_t$.


\section{Derivation of the first injectivity radius estimate}
\label{proofinject}

In this section we provide a proof of Theorem~\ref{inject}.

\

\noindent{\sl Step 1. Radius of definition of the exponential map. } 
First of all, we note that the injectivity radius of the Riemannian metric $g|_{\Sigma_0}$ induced on the initial
hypersurface $\Sigma_0 = t^{-1}(0)$ is controled, as follows.
Using Assumptions $(A3)$ and $(A4)$, we see that the Riemann curvature of the metric 
$g|_{\Sigma_0}$ is bounded and the volume of the unit geodesic ball $\vol_{g|_{\Sigma_0}}(\Bcal_{\Sigma_0}(p,1))$
is bounded below.
Therefore, according to \cite{CGT}, there exists a constant $i_1= i_1(K_2,v_0)$
such that the injectivity radius of $g|_{\Sigma_0}$ at the point $p$ is $i_1$ at least:
$$
\Inj_{g|_{\Sigma_0}}(\Sigma_0, p) \geq i_1.
$$
Moreover, according to \cite{JostKarcher} we can also assume that $i_1$ is sufficiently small so that, 
given any $\eps>0$  
there exists a coordinates $(x^\alpha)$ defined in a ball with definite size near $p$,
with $x^\alpha(p)=0$, such that the metric $g|_{\Sigma_0}$ is close to the
$n$-dimensional Euclidian metric $g_{E'}=\delta_{ij}$ (in these coordinates). 
More precisely, on the initial slice $\Sigma_0$ we have
$$
e^{-\eps} \, \delta_{ij} \leq g_{ij}(0,x^1, \ldots, x^n) \leq e^{\eps}  \, \delta_{ij},
\quad (x^1, \ldots, x^n) \in B_{E'}(0, i_1),
$$
where we have set $B_{E'}(0,r) := \big\{(x^1)^2 + \ldots + (x^n)^2 < r^2 \big\} \subset \R^n$. 
The latter can be regarded as a subset of $\Sigma_0$ by identifying a point with its coordinates, 
and we also use the notation $\Bcal_{E'}(p,r)$ for this Euclidian ball.

We can next introduce some coordinates $(x^\alpha) = (t,x^j)$ on the manifold, by propagating  
the coordinates $(x^j)$ chosen on $\Sigma_0$ along the integral curve of the vector field $T$. 
This construction allows us to cover the domain $\Omega$. 
From Assumption $(A2)$ (together with $(A1)$ and \eqref{Lie}) we deduce that the induced metric
on {\sl each slice} of the foliation is comparable with the $n$-dimensional Euclidian metric
in some time interval $[-i_2,i_2]$, that is 
$$
\aligned 
& (e^{-\eps} - K \, i_2) \, \delta_{ij} \leq g_{ij}(x) \leq (e^\eps + K \, i_2) \, \delta_{ij},
\\
& x= (t, x^1, \ldots, x^n) \in [-i_2,i_2] \times B_{E'}(0, i_1),
\endaligned 
$$
for some $K>0$ depending only on $K_0,K_1,K_2$. 

We then restrict attention to a smaller radius $i_2=i_2(K_0,K_1, K_2) \leq i_1$ such that 
$e^{-\eps} - K \, i_2 >0$, 
and we pick up $c_1 \geq 0$ sufficiently large 
so that $e^{-c_1} \leq e^{-\eps} - K \, i_2 \leq e^\eps + K \, i_2 \leq e^{c_1}$.
In turn, in view of Assumption $(A1)$ on the lapse function $n$
and of the expression \eqref{refer} of the reference Riemannian metric $g_T$,
the above inequalities imply that the reference Riemannian metric
$g_T$ is comparable to the $(n+1)$-dimensional Euclidean metric:
$$
e^{-c_2} \, \delta_{\alpha\beta} \leq g_{T,\alpha\beta} \leq e^{c_2} \, \delta_{\alpha\beta},
\quad
x= (t, x^1, \ldots, x^n) \in [-i_2,i_2] \times B_{E'}(0, i_2)
$$
for some constant $c_2 \geq c_1$ depending upon $c_1$ and $K_0$. 

Introducing on the manifold the $(n+1)$-dimensional Euclidian metric $E$ (which we define
in the constructed coordinates $(x^\alpha)$ and is, of course, independent of the point 
on the manifold)
and the corresponding Euclidian metric ball $\Bcal_E(p,i_2)$, we have established
\be
\label{equiva}
e^{-c_2} \, g_E \leq g_{T,q} \leq e^{c_2} \, g_E,
\quad
q \in \Bcal_E(p, i_2).
\ee
In the following we use the notation $|X|_E$ for the Euclidian norm of a vector $X$.

\

Our first task is to determine the radius of a ball on which the exponential map is well-defined.
This radius depends upon the reference vector field $T$. 
Let $\gamma:[0,s_0] \to M$ be a geodesic associated with the Lorentzian metric $g$
and  satisfying $\gamma(0)=p$. Assume that this geodesic is included in the Euclidian ball $\Bcal_E(p,i_2)$
(in which we have a good control of the metric $g_T$).
Obviously, we have
$$
\langle\gamma'(s), \gamma'(s) \rangle_g = \langle\gamma'(0), \gamma'(0) \rangle_g,
\qquad
s\in [0,s_0].
$$
On the other hand, to determine the length of $\gamma'(s)$ with respect to the {\sl reference} metric $g_T$,
we proceed as follows:
\begin{equation*}
\begin{split}
\big| \frac{d}{ds}\langle\gamma'(s), \gamma'(s)\rangle_T \big|
=
\big| \nabla_{T,\gamma'(s)} \big( g_T(\gamma'(s), \gamma'(s))\big) \big|
& = 2 \, \big| \langle \nabla_{g_T,\gamma'(s)}\gamma'(s), \gamma'(s)\rangle_T \big|
\\
& = 2 \, \big| \langle(\nabla_{g_T} - \nabla_g)_{\gamma'(s)} \, \gamma'(s), \gamma'(s)\rangle_T \big|
\\
& \leq 2 \, | \nabla_{g_T} - \nabla_g |_T \, |\gamma'(s)|_T^3.
\end{split}
\end{equation*}
So, by Lemma~\ref{Gamma}, $\big| \frac{d}{ds} |\gamma'(s)|_T^2 \big| \leq 2 \, K_3 \, |\gamma'(s)|_T^3$,
and, in consequence,
$$
\big| \frac{d}{ds} |\gamma'(s)|_T^{-1} \big|\leq K_3.
$$

By integration of the above inequality and {\sl provided} $s$ is small enough so that 
$2 s \, K_3 \, |\gamma'(0)|_T < 1$, we see that
\be
\label{geo1}
\frac{1}{2} \, |\gamma'(0)|_T \leq|\gamma'(s)|_T \leq 2 \, |\gamma'(0)|_T.
\ee
In view of \eqref{equiva} this implies
\be
\label{geo2}
\frac{e^{-c_2}}{2} \, |\gamma'(0)|_E \leq |\gamma'(s)|_E \leq 2 \, e^{c_2} \, |\gamma'(0)|_E.
\ee
These inequalities hold for all $s\in [0,1 / (2 K_3 \, |\gamma'(0)|_T)]$ as long as
$\gamma(s)\in \Bcal_E(p,i_2)$.
In particular, by restricting attention to geodesics whose initial vector has unit Euclidian length,
$|\gamma'(0)|_E=1$, we see that $\gamma([0,r_2]) \subset \Bcal_E(p,i_2)$
where $r_2 := i_2 e^{-c_2}/2$.
In turn, this establishes
that the exponential map  at the point $p$ is well-defined on the ball $B_E(0,r_2)$ 
with a range included in the geodesic ball $\Bcal_E(p,i_2)$.

\

\noindent{\sl Step 2. Conjugate radius estimate.}
Our second task is to determine a ball on which the exponential map is a local diffeomorphism, and 
we therefore need to control the length of a Jacobi field along a geodesic. 
Let $\gamma:[0, r_2] \rightarrow M$ be a $g$-geodesic satisfying $\gamma(0)=p$ and $|\gamma'(0)|_E = 1$. 
By the discussion in Step~1 we already know that the curve $\gamma$ lies in $\Bcal_E(p,i_2)$ 
and that $\max_{s \in [0,r_2]} |\gamma'(s)|_T \leq 2 \, e^{2c_2}$.
Given an arbitrary Jacobi field along $\gamma$,  $J=J(s)$, satisfying
$$
\aligned
& J''(s) = - \Riem (J(s), \gamma'(s)) \gamma'(s),
\\
& J(0)=0, \qquad | J'(0) |_T = 1,
\endaligned
$$
we need to control its Riemannian length $F(s) := |J|_T(s)$, as stated in \eqref{Jacobi} below.

Let $[0,s_0]$ be the largest subinterval of $[0,r_2/4]$ in which the inequality $| J |_T \leq 1$ holds.
Using the equation satisfied by the Jacobi field and taking into account the curvature bound $(A3)$, 
we deduce that, in the interval $[0,s_0]$, 
\begin{equation*}
\begin{split}
\big| \frac{d}{ds}\langle \nabla_{\gamma'} J, \nabla_{\gamma'}J\rangle_T \big|
& = 2 \, \big| \langle\nabla_{g_T,\gamma'}\nabla_{\gamma'} J,\nabla_{\gamma'} J\rangle_T \big|
\\
& \leq 2 \, | \nabla_{g_T} - \nabla_g |_T |\gamma'|_T \, |\nabla_{\gamma'}J|_T^2
   + 2 \, K_2 \, |\gamma'|_T^2| J|_T \, |\nabla_{\gamma'}J|_T. 
\end{split}
\end{equation*}
With \eqref{geo1} and the covariant derivative estimate in Lemma~\ref{Gamma},
we obtain 
\be
\label{ineq}
\big| \frac{d}{ds} |\nabla_{\gamma'}J|_T \big|
\leq 4 \, K_3 \, |\nabla_{\gamma'}J|_T + 8 \, K_2.
\ee

We can next integrate \eqref{ineq} over an arbitrary interval $[0,s]\subset [0,s_0]$, 
use the initial condition on the Jacobi field, and obtain  
$$
1 + {2 K_2 \over K_3} \, (1 - e^{-4K_3s}) \leq |\nabla_{\gamma'} J |_T 
\leq 1 + {2 K_2 \over K_3} \, (e^{4K_3s} - 1).  
$$
Assuming that $r_2$ is small enough so that 
${2 K_2 \over K_3} \, (1 - e^{-4K_3s}) \leq 1/2$ 
and $ {2 K_2 \over K_3} \, (e^{4K_3s} - 1) \leq$   
we infer that 
\be
\label{Jac1}
{1 \over 2} \leq |\nabla_{\gamma'} J |_T \leq 2. 
\ee 
Hence, using this inequality and Lemma~\ref{Gamma} we find $\frac{d}{ds}| J|_T \leq 2 + 2 K_3\leq 1$, 
so that 
\be
\label{Jac7}
F(s) = |J|_T (s)\leq (2+2 K_3) \, s \leq (2+2 K_3) \, r_2. 
\ee
Further assuming that $(2+2 K_3) \, r_2\leq 1$ we conclude that $s_0 = r_2$.

Next, we want to improve the rough estimate \eqref{Jac7}.
Since
\begin{equation*}
\begin{split}
\frac{d}{ds}\langle \nabla_{\gamma'}
J,J\rangle_T
& =\langle\nabla_{g_T,\gamma'}\nabla_{\gamma'} J, J\rangle_T
+\langle\nabla_{g_T,\gamma'} J,\nabla_{\gamma'} J\rangle_T
 \end{split}
\end{equation*}
then by substituting the previous estimates of $|J|_T (s)$ and $|\nabla_{\gamma'}J|_T(s)$
and performing similar calculations as above, we get
\begin{equation*}
\begin{split}
e^{- c_3} \leq \frac{d}{ds}\langle \nabla_{\gamma'}
J,J\rangle_T \leq e^{c_3} 
\end{split}
\end{equation*}
for some constant $c_3>0$. By integration this implies
$$
e^{-c_3} \, s
\leq \langle \nabla_{\gamma'} J,J\rangle_T
\leq e^{c_3} \, s
$$
and we arrive at the following lower bound for the norm of the Jacobi field:
$$
\aligned
F(s)
& \geq {\big |\langle \nabla_{\gamma'} J,J\rangle_T \big| \over | \nabla_{\gamma'} J |_T}
\geq \frac{e^{-c_3} \, s}{2} 
\geq e^{-c_4} \, s
\endaligned
$$
for some $c_4>0$. 

On the other hand, using again the above estimates we have
\begin{equation*}
\begin{split}
\frac{d}{ds} F & \leq \frac{1}{F} (\langle \nabla_{g_T,\gamma'} J,J\rangle_T + K_3 \, F^2)
\\
& \leq {e^{c_4} \over s} \Big( e^{c_3} \, s + K_3 (2 + 2K_3)^2 \, s^2 \Big) 
 \leq e^{c_5} 
\end{split}
\end{equation*}
for some constant $c_5>0$. This leads to the upper bound
$$
F(s) \leq e^{c_5} \, s.   
$$

In summary, we have established that the norm of the Jacobi field is comparable with $s$:
\be
\label{Jacobi}
e^{-c_4} \, s \leq F(s)\leq e^{c_5} \, s,  
\qquad
s\in[0,r_2].
\ee
By the definition of Jacobi fields these inequalities are equivalent to controling
the differential of the exponential map, that is for $s\in[0,r_2]$
$$
e^{-c_4} \, |W|_{T}
\leq | d\expb_{p,s\gamma'(0)}(W)|_T \leq e^{c_5} \, |W|_{T}.
$$ 
We conclude that the pull back of the reference metric to the tangent space at $p$ satisfies
\be
\label{expo}
\aligned 
& e^{-c_4} \, g_{T,p} \leq \big(\expb_p\big)^\star g_T \leq e^{c_5} \, g_{T,p}
\\
& \text{ in the ball } B_T(0,r_2)\subset T_p M.
\endaligned 
\ee
In particular, since the conjugate radius of the Lorentzian metric is precisely defined
from the reference Riemannian metric,
these inequalities show that the conjugate radius of the exponential map is $r_2$ at least.

\

\noindent {\sl Step~3. Injectivity radius estimate. } We are now in a position to establish
that $ \Inj_g(M,p,T)\geq r_3 := r_2 e^{-c_2}/4$.
We argue by contradiction and assume that $\gamma_1 : [0,s_1]\rightarrow M$
and $\gamma_2 :[0,s_2]\rightarrow M$ are two distinct $g$-geodesics satisfying $\max (s_1,s_2) \leq r_3$
and
$$
\aligned
& \gamma_1(0)=\gamma_2(0)=p,
\quad
|\gamma'_1(0)|_T =|\gamma'_2(0)|_T = 1,
\\
& \gamma_1(s_1)=\gamma_2(s_2) =: q.
\endaligned
$$
We will reach a contradiction and this will establish that the injectivity radius 
is greater or equal to $r_3$ 
(as can be checked by using the fact that the exponential map is at least a local diffeomorphism). 

By Step~1 we know that $\gamma_1, \gamma_2\subset \Bcal_E(p,2 e^{2c_2} r_3)$. By concatenating these two curves,
we construct a geodesic loop containing $p$,
$$
\gamma=\gamma_2^{-1}\cup\gamma_1:[0,s_1+s_2]\rightarrow \Bcal_E(p,2 e^{2c_2} r_3),
$$
which need not be smooth at $p$ or $q$.
Since $\gamma$ is contained in the image of the ball $B_T(p,r_2)$ under the exponential map, 
we can define an homotopy of $\gamma$ with the origin ($x=0$), 
by setting (in the coordinates constructed earlier)
$$
\Gamma_\eps(s)=\eps \gamma(s), \quad \eps \in [0,1].
$$
The curves $\Gamma_\eps:[0,s_1+s_2]\rightarrow \Bcal_E(p,2 e^{2c_2} r_3)$ satisfy
$$
\Gamma_\eps(0) = \Gamma_\eps(s_1+s_2)=p, \quad \Gamma_0([0,1])=p,
\quad \Gamma_1=\gamma.
$$

Moreover, we have $|\Gamma'_\eps(s)|_E \leq \eps 2 e^{2c_2} \leq 2 e^{2c_2}$ 
and thus $|\Gamma'_\eps(s)|_T \leq 2 e^{3c_2}$. In particular, the
$g_T$-lengths (computed with the reference metric) of the loops $\Gamma_\eps$ are less than
$$
L(\Gamma_\eps, g_T) \leq 2 e^{3c_2} r_3 = {r_2 \over 2}.
$$
 
Since the exponential map is a local diffeomorphism from the ball
$B_T(0,r_2)\subset T_pM$ to the manifold, and in view of the
estimate \eqref{expo} on the exponential map, it follows that all
the loops $\Gamma_\eps$ can be lifted to the ball
$B_T(0,r_2)$ in the tangent space with the {\sl same} origin
$0$. Consequently, we obtain a {\sl continuous} family of curves
$\Gammat_\eps:[0,s_1+s_2]\rightarrow T_pM$ satisfying
$$
\Gammat_\eps(0)=0, \quad \eps\in[0,1].
$$
At this juncture we observe that, since $\Gammat_\eps(s_1+s_2)$
(for $\eps\in[0,1]$) all cover the same point $p$ and since the
curve $\Gammat_0$ is trivial and the family is continuous,
$$
\Gammat_\eps(s_1+s_2)=0, \quad \eps\in[0,1].
$$
It remains to consider the lift of the original geodesic loop $\gamma$: under the lifting
the geodesics $\gamma_1, \gamma_2$ are sent to two distinct {\sl line segments}
(with respect to the vector space structure) originating at the origin $0$ which obviously do not
intersect. This is a contradiction
and we conclude that, in fact,  $\Inj_g(M,p,T)\geq r_3$ as announced.
This completes the proof of Theorem~\ref{inject}.


\section{Convex functions and convex neighborhoods}
\label{convexsection}

We establish now the existence of convex functions and convex neighborhoods in $M$.
Let us recall first some basic definitions. A function $u$ is said to be {\sl geodesically convex}
if the composition of $u$ with any geodesic is a convex function (of one variable).
A set $\Omega'\subset \Omega''$ is said to be {\sl relatively geodesically convex in $\Omega''$}
if, given any points $p,q \in \Omega'$ and any geodesic (segment)
$\gamma$ from $p$ to $q$ contained in $\Omega''$, one has $\gamma \subset \Omega'$.
A set $\Omega'$ is said to be {\sl geodesically convex in $\Omega''$} if $\Omega'$ is relatively
geodesically convex in
$\Omega''$ and, in addition, for any $p,q$, there exists a unique geodesic
$\gamma$ connecting $p$ and $q$ and lying in $\Omega'$.

We denote by $d_T$ the distance function associated with the reference Riemannian metric $g_T$.

\begin{theorem}[Existence of geodesically convex functions]
\label{convex}
Let $(M,g)$ be a differentiable $(n+1)$-manifold endowed with a Lorentzian metric 
$g$, satisfying the regularity assumptions $(A1)$-$(A4)$ for some point $p\in M$
and some future-oriented, unit, time-like vector field $T$, 
and let $g_T$ be the reference Riemannian metric associated with 
Then, for any $\eps \in (0,1)$ there exists a positive constant $r_0$ depending only upon $\eps$,
the foliation bounds $K_0,K_1$, the curvature bound $K_2$, the volume bound $v_0$,
and the dimension of the manifold and there exists a smooth function $u$ defined on $\Bcal_T(p,r_0)$
such that 
$$
\aligned
& (1-\eps) \, d_T(p,\cdot)^2 \leq u \leq (1+\eps) \, d_T(p,\cdot)^2,
\\
& (2-\eps) \, g_T \leq \nabla^2 u \leq (2+\eps) \, g_T.
\endaligned
$$ 
\end{theorem}

Hence, the function $u$ above is equivalent to the Riemannian distance function from $p$
and is geodesically convex for the Lorentzian metric. In the proof given below, the function $u$ 
is the Riemannian distance function associated with a new
 time-like vector field (denoted by $N$ in the proof
below).
The following corollary is immediate and provides us with a control of the radius of convexity,
which generalizes Whitehead theorem from Riemannian geometry \cite{Whitehead,CheegerEbin}.

\begin{corollary}[Existence of geodesically convex neighborhoods]
Under the assumptions of Theorem~\ref{convex}, for any $0<r<r_0$ there exists a set $\Omega_r \subset \Omega$
which is geodesically convex in $\Bcal_T(p,2r_0)$ and satisfies
$$
\expb_p(B_T(0,r))\subset \Omega_r \subset \expb_p(B_T(0,(1+\delta)r)).
$$ 
Moreover, one can always choose $\Omega_r$ so that
$$
\Bcal_T(p,r)\subset \Omega_r \subset \Bcal_T(p,(1+\delta)r),
$$
where $\Bcal_T(p,r)$ is the geodesic ball determined by the reference Riemannian metric.
\end{corollary}

\begin{proof}[Proof of Theorem~\ref{convex}]
 {\sl Step 1. Synchronous coordinate system. }
Given $\eps>0$, by applying the injectivity radius estimate in Theorem~\ref{inject}
to points near $p$, we see that there exists a constant $r_0$ depending on $K_0,K_1,K_2, v_0,\eps, n$
such that for any $q \in \Bcal_T(p,2r_0)$ the injectivity radius at $q$ is $2r_0$ at least,
and we can assume that
$$
e^{-\eps}\, g_{T,q} \leq (\expb_q)^\star g_T \leq e^\eps \, g_{T,q}, \qquad B_T(0,r_0) \subset T_q M, 
\, \, 
q \in \Bcal_T(p, 2r_0). 
$$ 

Let $\gamma=\gamma(s)$ be the (backward) time-like geodesic satisfying
$\gamma(0)=p$ and $\gamma'(0)=-T_p$, and consider the (past) point $q := \gamma(r_0/2)$.
The future null cone at $q$ with radius $r_0$  (the orientation being determined by the vector field $T$) is defined by
$$
C_q(r_0) := \big\{ V\in T_qM \, \big/ \,  |V|_{g_{T,q}} < r_0, \, |V|^2_{g_q} < 0, \, \langle V,T\rangle>0 \big\}.
$$
Observe that the $g_T$-length of $\gamma$ between $p$ and $q$ is approximatively
$r_0/2$ and that the norm $|\gamma'|_T$ is almost $1$, while
$|\gamma'(q)|_{g_q}^2=1$ and $\langle -\gamma',T\rangle_g>0$.
By the injectivity radius estimate in Theorem~\ref{inject}
the exponential map at $q$ is a diffeomorphism from $C_q(r_0)$ onto its image which, moreover, contains the
original point $p$.

Next, introduce the set of vectors that are ``almost'' parallel to $T$:
$$
C_q(r_0,\eps) := \big\{ V \in T_qM \, \big/ \,  |V|_{T,q} < r_0,  \,
\langle V, T\rangle_{g_q} > 0,
\, \frac{\langle V,V\rangle_{g_q}}{\langle V, V \rangle_{T,q}} > 1 - \eps \big\}.
$$
The notation $c(\eps)>0$ is used for constants that depend only on $K_0,K_1, K_2, v_0, n,\eps$
and satisfy $\lim_{\eps \to 0} c(\eps)=0$. We claim that there is constant $c(\eps)>0$ such that
\be
\label{58}
\Bcal_T(p,c(\eps) r_0) \subset \expb_q(C_q(r_0,\eps)).
\ee
Actually, we have $\Bcal_T(p,c(\eps)r_0)\subset \Bcal_T(q,(\frac{1 }{2}+c(\eps))r_0)$, hence
$$
\Bcal_T (p,c(\eps)r_0) \subset \expb_q \big( B_T(0,(\frac{1}{2}+c(\eps))r_0) \big).
$$
Since the metrics $g_{T,0}$ and $g_{T,q}$ are comparable (under the exponential map at $q$)
we see that geodesics $\sigma$ connecting $q$ and points of $\Bcal_T(p,c(\eps)r_0)$ make an angle
$\leq c(\eps)$ with $-\gamma'(q)$ at the point $q$ (as measured by the metric $g_{T,q}$).
By reducing the constant $c(\eps)$ if necessary, the claim is proved.
 
Let $\tau$ be the Lorentzian distance from $q$: it is defined on $\expb_q(C_q(r_0))$ and is
 a smooth function on $\expb_q(C_q(r_0))\setminus \{p\}$.
Using the claim \eqref{58} we deduce that $\tau$ is smooth in the ball $\Bcal_T (p,c(\eps){r_0})$ and
satisfies
\be
\label{tau}
\big( \frac{1}{2} - c(\eps) \big) \, r_0 < \tau
< \big(\frac{1}{2} + c(\eps) \big) \, r_0
\quad \text{ in the ball } \Bcal_T (p,c(\eps){r_0}).
\ee
It is clear also that
$$
|\nabla \tau|_{g}^2=-1, \qquad \nabla^2\tau (\nabla \tau,\cdot)=0.
$$

We now introduce a new foliation. 
Let $(z^j)$ be coordinates on the level set hypersurface $\tau=\tau(p)$, and
by following the integral curves of the (unit, time-like) vector field
$$
N := \nabla \tau
$$
let us construct coordinates $(z^\alpha)$ with $z_0:=\tau$
in which the Lorentzian metric $g$ takes the simple form
$$
g = - (dz^0)^2 + g_{ij} \, dz^idz^j.
$$
Let $g_N = \langle \cdot, \cdot\rangle_N$ be a (new) reference Riemannian metric based on the vector field $N$.

By Lemma~\ref{Gamma} 
using the equation satisfied by (future) $g$-geodesics $\sigma$ we see that
$$
\big| \frac{d}{d\tau}\log |\sigma'(\tau)| \big| \leq K_3 \, r_0.
$$
(Recall that we allow $r_0$ to depend upon $\eps$.)
This inequality shows that the vector field $N$ makes an angle $\leq c(\eps)$ with $T$,
everywhere on $\expb_q(C_q(r_0,\eps))$. From this, we conclude that the two metrics are comparable:
$$
(1-c(\eps)) \, g_T \leq g_N \leq (1+c(\eps)) \, g_T  \quad \text{ in the cone } \expb_q(C_q(r_0,\eps)).
$$

\

\noindent{\sl Step 2. Hessian comparison theorem and curvature bound for the reference metric $g_N$.}
Since $p\in \expb_q(C_q(r_0))$, let $\sigma:[0,\tau(p)] \to M$ be the future time-like
geodesic connecting $q$ to $p$, and let $\Vt$ be the Jacobi field defined along $\sigma$ such that
$$
\Vt(0)=0, \quad \Vt(\tau(p))=V, 
$$
where $V\in T_pM$ satisfies the orthogonality condition $\langle\nabla \tau,V\rangle=0$. 
Then we have
\begin{equation*}
\begin{split}
-\nabla^2\tau(V,V)&=-\langle\Vt,\nabla_{\nabla
\tau}\Vt\rangle=\langle\Vt,\nabla_{\frac{\del}{\del
\tau}}\Vt\rangle\\
&=\int_0^{\tau(p)}\langle\nabla_{\frac{\del}{\del
\tau}}\Vt,\nabla_{\frac{\del}{\del
\tau}}\Vt\rangle_{g}-R(\sigma',\Vt,\sigma',\Vt) =: I(\Vt,\Vt).
\end{split}
\end{equation*}
Recall that in the absence of conjugate points Jacobi fields minimize the index form $I(V,V)$
among all vector fields with fixed boundary values.
By applying a standard comparison technique from Riemannian geometry
on the orthogonal space $(\nabla \tau)^\perp$
(on which the Lorentzian metric induces a Riemaniann metric)
we control the Hessian of $\tau$ in terms of the curvature bound $K_2$:
\be
\label{Hess}
\frac{\sqrt{K_2(1+c(\eps))}}{\tan \sqrt{K_2(1+c(\eps))}\tau} \, g\mid_{(\nabla \tau)^\perp} 
\leq (-\nabla^2\tau)|_{(\nabla
\tau)^{\perp}}\leq \frac{\sqrt{K_2(1+c(\eps))}}{\tanh \sqrt{K_2(1+c(\eps))}\tau} \, g\mid_{(\nabla \tau)^\perp}.
\ee
Since $-\nabla^2_{ij}\tau=\frac{1}{2}\frac{\del g_{ij}}{\del \tau}$, we deduce from \eqref{Hess}
that
\be
\label{40}
\frac{g_{ij}}{\tau}\leq \frac{\del g_{ij}}{\del \tau}\leq
\frac{3g_{ij}}{\tau} \quad \text{ in the cone } \expb_q(C_q(r_0)).
\ee

Combining \eqref{40} with the curvature formulas derived in Section~3, i.e.
\begin{equation*}
\begin{split}
R_{ijkl}&=R^{\Sigma}_{ijkl}+\frac{1}{4} \Big( \frac{\del g_{ik}}{\del \tau}\frac{\del g_{jl}}{\del \tau}
     -\frac{\del g_{il}}{\del \tau}\frac{\del g_{jk}}{\del \tau}\Big),
\\
R_{0jil}&={1 \over 2} \, \big( \nabla_{l}(\frac{\del}{\del \tau}g_{ij})-\nabla_i(\frac{\del}{\del t}g_{lj}) \big),
\\
R_{i0j0} &=-{1 \over 2}\frac{\del^2 g_{ij}}{\del \tau^2}
+\frac{1}{4}g^{pq}\frac{\del g_{ip}}{\del \tau}\frac{\del g_{jq}}{\del \tau},
\end{split}
\end{equation*}
we conclude that
\be
\label{timederiv}
\big| \frac{\del^2 g_{ij}}{\del \tau^2} \big| \leq \frac{C}{\tau^2} \quad \text{ on } \expb_q(C_q(r_0)).
\ee

Finally, relying on the formulas for the curvature of the reference Riemannian metric $g_N$, we obtain
$$
|\Riem_{g_N}|_{N}\leq \frac{C}{\tau^2} \quad \text{ on } \expb_q(C_q(r_0)).
$$
(Observe that, as could have been expected, the upper bound blows-up as one approach the point $q$
which is the base point in our definition of the distance.) 
In particular, this implies the following curvature bound near the point $p$:
$$
|\Riem_{g_N}|_{N}\leq Cr_0^{-2}  \quad \text{ on the ball } \Bcal_T(p,c(\eps)r_0).
$$

\

\noindent {\sl Step 3. Constructing geodesically convex functions.}
Since  the metrics $g_T$ and $g_N$ are comparable, the volume ratio
$(1/r_0^{n+1}) \, \vol_{g_N} \Bcal_N(p,c(\eps)r_0)$ is uniformly bounded (above and) below.
By applying the theory for Riemannian metrics in \cite{CGT}, the injectivity radius of the metric $g_N$
is bounded from below by $c(\eps)r_0$. Let
$$
u(x) := d_{g_N}(p,x)^2
$$
be the (square) of the distance function associated with the Riemannian metric $g_N$,
which is a smooth function defined on the geodesic ball $\Bcal_N(p,c(\eps)r_0)$.
By the standard Hessian comparison theorem for Riemannian manifold we have
$$
(2-\eps) \, g_{N,\alpha\beta} \leq \nabla_{g_N,\alpha} \nabla_{g_N,\beta} u
\leq
(2+\eps) \, g_{N,\alpha\beta} \qquad \text{ on the ball } \Bcal_N(p,c(\eps)r_0).
$$

In terms of the original Lorentzian metric $g$, the Hessian of the function $u$ is
$$
\nabla_\alpha \nabla_\beta u = \nabla_{g_N,\alpha}\nabla_{g_N,\beta} u
        + (\Gamma^\gamma_{g_N,\alpha\beta}-\Gamma^\gamma_{\alpha\beta}) \,\frac{\del u}{\del x^\alpha}.
$$
Since
$| \Gamma_{g_N} - \Gamma|_{N} \leq C \, \sup |\frac{\del g_{ij}}{\del \tau}|\leq C'$ by the estimate \eqref{40}
and since also $|\nabla u|_{N}\leq 2 \, d_{g_N}$ on $\Bcal_N(p,r_0)$, we conclude that
$$
(2- \eps ) \, g_{N,\alpha\beta}
\geq {\nabla}_\alpha{\nabla}_\beta u \geq (2+\eps) \, g_{N,\alpha\beta}
\qquad \text{ in the ball } \Bcal_{N}(p,c(\eps)r_0).
$$
This completes the proof of Theorem~\ref{convex}.
\end{proof}


\section{Injectivity radius of null cones}
\label{nullsection}

We turn our attention now to null cones within foliated Lorentzian manifolds.
Our main result (Theorem~\ref{nulltheorem} below) provides a lower bound for the null injectivity 
radius under the main assumption that the exponential map is defined in some ball and 
the null conjugate radius is already controled. Hence, contrary to the presentation in 
Section~\ref{WS-0} our main assumption (see $(A3')$ below) is not directly stated as a curvature bound.
However, under additional assumptions, 
it is known that the conjugate radius estimate can be deduced from an $\Lbf^p$ curvature bound, 
so that our result is entirely relevant for the applications. 

Indeed, in a series of fundamental papers \cite{KR1,KR2,KR4}, Klainerman and Rodnianski
assumed on an $\Lbf^2$ curvature bound and estimated the null conjugate and injectivity radii
for Ricci-flat Lorentzian ${(3+1)}$-manifolds. Our result in the present section is a continuation of 
the recent work \cite{KR4} and covers a general class of 
Lorentzian manifolds with arbitrary dimension, while our proof is local and geometric and so conceptually simple.

We use the terminology and notation introduced in Section~\ref{LO-0}. In particular,
a point $p\in M$ and a reference vector field $T$ are given, and
$N_p^-$ denotes the past null cone in the tangent cone at $p$.
The null exponential map $\expb^N_p : B_T^N(0,r) \to M$ is defined over a subset of this cone,
$$
B_T^N(0,r) := B_T(0,r)\cap N_p^-,
$$
and allows us to introduce the (past) null injectivity radius $\NullInj_g(M,p,T)$. We also set
$$
\Bcal_T^N(p, r) := \expb_p^N(B_T^N(0,r)).
$$

We consider a domain $\Omega \subset M$ containing some point $p$
on a final slice $\Sigma_0$ and foliated as
\be
\label{foliation11}
\Omega = \bigcup_{t \in [-1,0]} \Sigma_t,
\qquad p \in \Sigma_0.
\ee
We assume that there exist positive constants $K_0, K_1, K_2$ such that
\begin{equation*}
\tag{A1} 
e^{-K_0} \leq n \leq e^{K_0} \quad \text{ in } \Omega, 
\end{equation*}
\begin{equation*}
\tag{A2}
|\Lie_Tg |_T \leq K_1 \quad \text{ in } \Omega,
\end{equation*}
the null conjugate radius at $p$ is $r_0$ (at least) and the null exponential map satisfies
\begin{equation*}
\tag{A3'}
e^{-K_2} \, g_{T,p} \mid_{B_T^N(0,r_0)} \leq
\big({\expb^N_p}\big)^\star ( g_T \mid_{\Bcal_T^N(0,r_0)})
\leq 
e^{K_2} \, g_{T,p}\mid_{B_T^N(0,r_0)}
\end{equation*}
and, finally, there exists a coordinate system on the initial slice $\Sigma_{-1}$
such that the metric $g\mid_{\Sigma_{-1}}$ is comparable to the $n$-dimensional Euclidian metric $g_{E'}$
in these coordinates:
\begin{equation*}
\tag{A4'}
e^{-K_0} \, g_{E'} \leq g\mid_{\Sigma_{-1}} \leq e^{K_0} \, g_{E'}
\quad \text{ in } \Bcal_{\Sigma_{-1}, E'}(p,r_0). 
\end{equation*}
We refer to $K_2$ as the {\sl effective conjugate radius constant.}

\begin{theorem}[Injectivity radius of null cones]
\label{nulltheorem} Let $M$ be a differentiable ${(n+1)}$-manifold,
endowed with a Lorentzian metric $g$ satisfying the regularity assumptions
$(A1)$, $(A2)$, $(A3')$, and $(A4')$ at some point $p$ and for some foliation \eqref{foliation1}.
Then, there exists a positive constant $i_0$ depending only upon the foliation bounds $K_0,K_1$,
the null conjugate radius $r_0$, the effective conjugate radius constant $K_2$, and the dimension $n$ such
that the null injectivity radius of the metric $g$ at $p$ satisfies
$$
\NullInj_g(M,p,T) \geq i_0.
$$
\end{theorem}

It is interesting to compare the assumptions above with the ones in Section~\ref{WS-0}.
Assumptions $(A1)$ and $(A2)$ are concerned with the property of the foliation and
were already required in Section~\ref{WS-0}.

Assumption $(A3')$ should be viewed as a weaker version of the $L^\infty$ curvature condition $(A3)$.
Recall that, under the assumptions of Theorem~\ref{inject} which included a curvature bound,
an analogue of $(A3')$ valid in the whole of $\Omega$ was already established in \eqref{expo}.
It is expected that $(A3')$ is still valid when the curvature in every spacelike slice
is solely bounded in some $L^m$ space.

Indeed, at least when the spatial dimension is $n=3$ and the manifold is Ricci-flat,
according to Klainerman and Rodnianski \cite{KR1,KR2} 
Assumption $(A3')$ is a consequence of the following $\Lbf^2$ curvature bound
\be
\label{Lsquare22}
\|\Riem_g\|_{\Lbf^2(\Sigma_{-1},g_T)} \leq K_2'
\ee
for some constant $K_2'>0$.

Assumption $(A4')$ concerns the metric on the initial hypersurface and is only slightly stronger
than the volume bound $(A4)$. Furthermore, according to 
Anderson \cite{Anderson0} and Petersen \cite{Petersen} 
the property $(A4')$ is also a consequence of the curvature bound  
\be
\label{Lm}
\|\Riem_g\|_{\Lbf^m(\Sigma_{-1},g_T)} \leq K_2'
\ee
for $m > n/2$ and some constant $K_2'>0$ and a volume lower bound at every scale
\be
\label{volumes}
r^{-n} \, \vol_{g|_{\Sigma_0}}(\Bcal_{\Sigma_0}(p,r))\geq v_0,
\qquad 
r \in (0,r_0]. 
\ee

In summary, by combining Theorem~\ref{nulltheorem} above with the results in \cite{KR4,JostKarcher} we conclude:

\begin{corollary}[Einstein field equations of general relativity]
\label{ricci}
Let $(M,g)$ be a Lorentzian ${(3+1)}$-manifold satisfying the vacuum Einstein equation
\be
\label{Einstein4}
\Ric_g = 0.
\ee
Suppose that near some point $p \in M$ there exists a foliation $\Omega$ of the form
\eqref{foliation11} satisfying Assumptions $(A1)$-$(A2)$ and
such that the $\Lbf^2$ curvature assumption \eqref{Lsquare22} holds on the initial spacelike hypersurface
$\Sigma_{-1}$. Then, there exists a positive constant $i_0$ depending only upon the foliation bounds
$K_0,K_1$ and the curvature bound $K_2'$ such that the null injectivity radius satisfies
$$
\NullInj_g(M,p,T) \geq i_0.
$$
\end{corollary}

\begin{proof}[Proof of Theorem~\ref{nulltheorem}]
\, {\sl Step 1. Localization of the past null cone $\Ncal^-(p)$ between two flat null cones. }
Assumption $(A3')$ provides us with a bound on the null conjugate radius,
we need to control the null cut locus radius. 
We proceed as in Section~\ref{proofinject} and introduce coordinates near the point $p$
such that $x^\alpha(p) = 0$. Precisely,
relying on Assumptions $(A1)$, $(A2)$, and $(A4')$, we determine the coordinates
$x=(x^\alpha)$ so that $x^0=t$ and the spatial coordinates $(x^j)$ are transported
(via the gradient of the function $t$) from the coordinates prescribed on the initial slice $\Sigma_{-1}$.
Then, the Lorentzian metric reads $g = -n^2 \, dt^2 + g_{ij} \, dx^idx^j$ and satisfies for some $C_0, C_1>0$
\be
\label{CC}
{1 \over C_0} \leq n^2 \leq C_0,
\qquad
{1 \over C_1} \, \delta_{ij} \leq g_{ij} \leq C_1 \, \delta_{ij}, 
\ee
for all $-r_0< t\leq 0$ and $(x^1)^2 + \ldots + (x^n)^2 \leq (r_0)^2$, and in these coordinates
the reference Riemannian metric $g_T$ is comparable to the ${(n+1)}$-dimensional Euclidian metric
$g_E : = dt^2 + (dx^1)^2 + \ldots (dx^n)^2$:
\be
\label{found}
{1 \over C_1} \, g_E \leq g_T \leq C_1 \, g_E.
\ee
Denote by $\Bcal_E(q,r)$ the Euclidean ball with center $q$ and radius $r$.
Note that these inequalities holds within a neighborhood of $p$ in $\Omega$.
The forthcoming bounds will hold in a neighborhood of the past null cone only.
To simplify the notation we set
$$
c_0 := {1 \over C_0}, \qquad c_1 := {1 \over C_1}.
$$

In each time slice of parameter value $t=a$ we introduce the $n$-dimensional Euclidian ball with radius $b$
$$
\Acal^a_{<b} := \big\{ t = a, \quad (x^1)^2 + \ldots + (x^n)^2 < b^2 \big\} \subset \Sigma_a,
$$
which is centered around the point $p'$ with coordinates $(a,0,\cdots,0)$.
We also define $\Acal^a_{>b}$, $\Acal^a_{[c,d]}$,\ldots similarly.

For any point $q$ in a slice $\Sigma_{t_0}$ satisfying $-r_0\leq t_0<0$ and $x^1(q)^2 + \cdots+x^n(q)^2 <
c_1^2 \, t_0^2$
we consider the line (for the Euclidian metric) connecting $q$ to $p$:
$$
\gamma(\tau) = \big( \tau,\frac{\tau}{t_0}x^1(q),\cdots,\frac{\tau}{t_0}x^n(q)\big),
\qquad \tau\in [t_0,0].
$$
This is a timelike curve for the Lorentzian metric $g$, since
\begin{equation*}
\begin{split}
|\gamma'(\tau)|^2= - n^2 + g_{ij} \, \frac{x^i(q)}{t_0} \, \frac{x^j(q)}{t_0}
 < - c_0 + c_1 < 0,
\end{split}
\end{equation*}
which shows that
$$
\Acal^t_{< c_1 |t|}\subset \Ical^-(p), \qquad  t \in (-r_0,0).
$$

On the other hand, we claim that the larger Euclidian cone $\Acal^t_{< C_1 \, |t|}$
contains the null cone, in other words 
$$
\Acal^t_{\geq C_1 \, |t|}\subset (\mathcal{N}^-(p)\cup \Ical^-(p))^{c},
\qquad t \in (- c_1 \, r_0,0).
$$
Indeed, arguing by contradiction we suppose there exist  a time $t_0 \in (- c_1 \, r_0,0)$
and  a point $q\in \Acal^{t_0}_{\geq C_1 \, t_0}$
connected to $p$ by a causal curve $\gamma=\gamma(s)$ with $\gamma(0)=p$.
After reparametrizing (in time) the curve is necessary we can assume that
$\gamma(\tau) = (\tau,x^j(\tau))$ for some $t_0' \leq \tau \leq 0$,
as long as the point $\gamma(\tau)$ lies in the coordinate system under consideration.
For this part of the curve at least we have
$$
0 \geq |\gamma'|^2=-n^2+g_{ij}\frac{dx^i}{d\tau}\frac{dx^j}{d\tau},
$$
which by \eqref{CC}
implies that $(\frac{dx^1}{d\tau})^2 + \ldots + (\frac{dx^n}{d\tau})^2 < C_1 \, C_0$. Therefore,
after integration we find
$$
\big( x^1(q)^2 + \cdots+x^n(q)^2\big)^{1/2} (t_0')
\leq
\sqrt{C_0 \, C_1} \, t_0'\leq \sqrt{C_0 c_1} \, r_0 < r_0.
$$
Hence, we can choose $t_0'=t_0$, the whole curve lies in the system of coordinates,
and is parametrized in the form $\gamma(\tau) = (\tau,x(\tau))$, $(\tau\in [t_0,0])$.
Moreover, we have $|x(t_0)|\leq \sqrt{C_1 \, C_0} \, |t_0| < C_1 \, |t_0|$,
which contradicts our assumption $q\in \Acal^{t_0}_{\geq C_1 \, t_0}$.

In conclusion, we have localized the slices of the past null cone within ``annulus'' regions:
$$
\Ncal^-(p) \cap \Sigma_t\subset \Acal^t_{[c_1 \, |t|, C_1|t|]},
\qquad t\in [-c_1 \, r_0,0].
$$

\

\noindent{\sl Step 2. The past null cone $\Ncal^-(p)$ viewed as a graph with bounded slope. }
We now obtain a Lipschitz continuous parametrization of the null cone.
For any fixed $q \in \Acal^{- c_1 r_0}_{\leq c_1^2 r_0}$ we consider the vertical curve
passing through $q$:
$$
\gamma_q(\tau)=(\tau,x^1(q),\cdots,x^n(q)), \qquad \tau\in [-c_1 r_0,0].
$$
By Step~1 we know that there exists $\tau_q$ such that $\gamma_q(\tau_q)\in \Ncal^-(p)$.
Moreover, $\tau_q$ is unique since $\Ncal^-(p)$ is achronal, and this defines a map
$$
F : \Acal^{-c_1 r_0}_{\leq c_1^2 r_0} \rightarrow \Ncal^-(p)
$$
such that $F(q) = \gamma_q(\tau_q)$. It is obvious $F(- c_1 r_0,0) = p$.

We claim that the map $F$ is Lipschitz continuous with Lipschitz constant less than $C_1$,
as computed with the Euclidean metric $E$. Namely, by contradiction, suppose that
$|F(q_1) - F(q_2)|_E > C_1 \, |q_1 - q_2|_E$ for some $q_1, q_2 \in \Acal^{- c_1 r_0}_{\leq c_1^2 r_0}$, then
by \eqref{found} in Step~1, $F(q_1)$ would be chronologically
related to $F(q_2)$ and this would contradict the fact that $\Ncal^-(p)$ is achronal.
Moreover, from Step~1 it follows that
$$
F(\Acal^{- c_1 r_0}_{\leq c_1^2 r_0}) \supset \Ncal^-(p) \cap \Bcal_E(p,c_1^3 r_0).
$$

\

\noindent{\sl Step 3. Constructing an homotopy of curves on the null cone. }
Suppose that $\gamma_1, \gamma_2$ are two (past) null geodesics from $p$ satisfying
$$
\aligned
& \gamma_1(0)=\gamma_2(0),
\qquad
|\gamma_1'(0)|_{T} = |\gamma_2'(0)|_{T} = 1,
\\
& \gamma_1(s_1)=\gamma_2(s_2).
\endaligned
$$
We claim that $\max(s_1, s_2) > c_1^6 r_0$, which will establish the desired injectivity bound
by setting $i_0 = c_1^6 r_0$.

We argue by contradiction and assume that $\max(s_1, s_2) < c_1^6 r_0$.
Taking into account Assumption~$(A2)$ and applying exactly the same arguments as
in Step~1 of Section~\ref{proofinject} we see that the $g_T$-lengths of the curves
$\gamma_1,\gamma_2$ satisfy
$$
L(\gamma_j,g_T) \leq s_j \, e^{C C_1 \, s_j} \leq c_1^{5+3/4} \, r_0 \qquad (j=1,2).
$$
By Step~1 of the present proof we know that the Euclidean lengths of $\gamma_1, \gamma_2$ satisfy
$$
L(\gamma_j,g_E) \leq c_1^{5+1/4} \, r_0 \qquad (j=1,2).
$$
In particular, $\gamma_1, \gamma_2 \subset \Ncal^-(p) \cap \Bcal_E(p, c_1^5 r_0)$ and we can thus
concatenate the curve $\gamma_1, \gamma_2$ and obtain 
$$
\gamma := \gamma_2^{-1} \cup \gamma_1: [0, s_1+s_2] \rightarrow
\Ncal^-(p) \cap \Bcal_E(p,c_1^5 r_0).
$$

Since $F(\Acal^{- c_1 r_0}_{\leq c_1^2 r_0})\supset \Ncal^-(p) \cap \Bcal_E(p, c_1^3 r_0)$,
there exists a smooth family of curves $\sigma_\eps:[0,s_1+s_2]\rightarrow \Ncal^-(p)$
such that
$$
\aligned
& \sigma_1 = \gamma, \quad \sigma_0=p,
\\
& \sigma_\eps(0) = \sigma_\eps(s_1 + s_2) = p, \qquad \eps\in [0,1].
\endaligned
$$
Specifically, we choose
$$
\sigma_\eps(s) := F( \eps F^{-1}\gamma(s)),
$$
where the multiplication by $\eps$ is defined by relying on the linear structure of
$\Acal^{- c_1 r_0}_{\leq c_1^2 r_0} \approx B_{\R^n}(0, c_1^2 r_0)$.
Equivalently, by setting $x^i(s)=x^i(\gamma(s))$ we have the explicit formula
$$
\sigma_\eps(s) = F\big(- c_1 r_0, \eps x^1(s), \cdots, \eps x^n(s)\big).
$$
It is clear that the Euclidean and $g_T$-lengths of $\sigma_\eps$ satisfy
$$
\aligned
& L(\sigma_\eps,g_E) \leq \eps (1 + C_1) \, L(\gamma,g_E) \leq c_1^{4+1/8} \, r_0,
\\
& L(\sigma_\eps,g_T) \leq c_1^{3+5/8} \, r_0.
\endaligned
$$

By Assumption $(A3')$ on the null conjugate radius, we can lift
to the null cone of the tangent space $T_pM$
the continuous family of loops $\sigma_\eps$, and we obtain
a continuous family of curves $\tsigma_\eps$ defined on $[0,s_1+s_2]$ such that
$$
\tsigma_\eps(0) = 0, \quad L(\tsigma_\eps,g_{T,p}) \leq c_1^3 \, r_0.
$$
Observe that the property $L(\tsigma_\eps,g_{T,p})\leq c_1^3 \, r_0 \preceq r_0$ guarantees
the existence of this continuous lift.
By continuity,  all of the curves $\tsigma_\eps$ are loops containing $0$.
As observed earlier in the proof for the case of bounded curvature,
$\tsigma_1$ consists of two distinct segments which clearly can not form a closed loop
and we have reached a contradiction.
\end{proof}


\section{Injectivity radius of an observer in a Lorentzian manifold}
\label{nopre}

\subsection*{Main result}

We are now a in a position to discuss and prove Theorem~\ref{nofoliation}
stated in the introduction.   
As we have seen in the proof of the previous section, once the injectivity radius is controled,
one can construct a foliation satisfying certain ``good'' properties. On the other hand, the concept
of injectivity radius is clearly independent of any prescribed foliation.
As this is more natural, we will now present a general result which avoids to assume a priori 
the existence of a foliation.
This will be achieved by relying on purely geometric and intrinsic quantities and
constructing coordinates adapted to the geometry. Such a result is conceptually very important
in the applications. The result and proof in this section should be viewed as a Lorentzian generalization
of Cheeger, Gromov, and Taylor's technique \cite{CGT}, originally developed for Riemannian manifolds.

Let $(M,g)$ be a differentiable ${(n+1)}$-manifold endowed with a Lorentzian metric tensor $g$, 
and consider a point $p\in M$ and a vector $T\in T_pM$ with $g_p(T,T)=-1$.
That is, we now fix a single observer located at the point $p$.  
As explained in Section~\ref{LO-0} 
the vector $T$ induces an inner product $g_T=\langle \, , \, \rangle_T$ on the tangent space $T_pM$.
We assume that the exponential map $\expb_p$ is defined in some ball $B_T(0,r_0)\subset T_pM$ 
determined by this inner product, which is of course always true in a sufficiently small ball. 
Controling the geometry at the point $p$ precisely amounts to estimating  
the size of this radius $r_0$ where the exponential map is defined and has some good property.  
We restrict attention to the geodesic ball $\Bcal_T(p,r_0) := \expb_p(B_T(0,r_0))$; 
recall that these sets depend upon the vector $T$ given at $p$.

As explained in the introduction, by $g$-parallel translating the vector $T$ at $p$ 
along a geodesic $\gamma$ from $p$, we can define get a future-oriented unit time-like vector field $T-\gamma$ 
defined along this geodesic. To this vector field and the Lorentzian metric $g$ we can 
associate a reference Riemannian metric $g_{T_\gamma}$ along the geodesic.
In turn, this allows us to compute the norm $| \Riem_g |_{T_\gamma}$ of the Riemann curvature tensor 
along the geodesic. 

Of course, whenever two such geodesics $\gamma, \gamma'$ meet away from $p$, 
the corresponding vectors $T_\gamma$ and $T_{\gamma'}$ are generally {\sl distinct.}  
If we consider the family of all such geodesics we therefore obtain a (generally) multi-valued vector field 
defined in the geodesic ball $\Bcal_T(p,r_0)$. We use the same letter $T$ to denote this vector field. 
In turn we can still compute the Riemann curvature norm $| \Riem_g |_T$ by taking into account every 
value of $T$.  

The key objective of the present section is the study of the geometry of the local covering
$\expb_p:B_T(0,r_0)\rightarrow \Bcal_T(p,r_0)$ and to compare the Lorentzian metric $g$ defined on the 
manifold $M$ with the reference Riemannian metrics $g_T$.
As we will see in the proof below, it will be convenient to pull the metric ``upstairs'' 
on the tangent space at $p$, using the exponential map. 
Indeed, this will be possible once we will have estimated the conjugate radius 
(in Step 1 of the proof below) and will know that the exponential map is non-degenerate on $B_T(0,r_0)$. 
Pulling back the Lorentzian metric $g$ on $M$ by the exponential map 
we get a Lorentzian metric $g=\expb_p^\star g$ defined in the tangent space, 
on the ball $B_T(0,r_0)$. We use the same letter $g$ to denote this metric. 
Then, the geometry in the tangent space is particularly simple, since 
the $g$-geodesics on $M$ passing through $p$ are radial straightline in $B_T(0,r_0)$. 

A third view point could be adopted by restricting attention within the cut-locus from the point $p$, 
and by imposing the curvature assumption within the cut-locus only. 

We are in a position to prove the main result of the present paper that was stated in Theorem~\ref{nofoliation}.

\begin{proof}[Proof of Theorem~\ref{nofoliation}]
After scaling we may assume that $r_0=1$, and so we need to show 
\be
\label{toprove}
\Inj_g(M,p,T)\geq c(n) \, \vol_{g}(\Bcal_T(p,c(n))).
\ee

\

\noindent{\sl Step 1. Estimates for the metric $g_T$ and its covariant derivative.} 
Let $E_0=T,$ $E_1,\cdots,E_n$ be an orthonormal frame at the origin in $T_pM$ for the Lorentzian metric $g$.
By $g$-parallel transporting this basis along along a radial geodesic 
$\gamma=\gamma(r)$, satisfying $\gamma(0)=0$, $|\gamma'(0)|_{T}=1$,
we get an orthonormal frame defined along the geodesic.  
We use the same letters $E_\alpha$ to denote these vector fields.
Since
$$
\frac{d}{dr}\langle E_\alpha, E_\beta \rangle_{g} = 0, 
$$
we infer that
$$
|E_i|_T^2 = |E_i|_g^2 = 1 \quad \text{ along the geodesic.} 
$$
The same argument also implies
\be
\label{gamm}
|\gamma'(r)|_{T}^2=|\gamma'(0)|_{T}^2=1, 
\qquad
|\gamma'(r)|_{g}^2=|\gamma'(0)|_{g}^2=1, 
\ee
and $\gamma'(r) = c^\alpha E_\alpha(r)$ with constant (in $r$) scalars $c^\alpha$ and 
$\sum |c^\alpha|^2 = |\gamma'(0)|_T = 1$. We used here that, by definition, 
$\gamma'$ is $g$-parallel transported. 

Let $V = a^\alpha(r) \, E_\alpha(r)$ be a Jacobi field along a radial geodesic $\gamma=\gamma(r)$,
with $V(0)=0$ and $|V'(0)|_T = 1$. Then, the Jacobi equation takes the form
$$
(a^\alpha)''(r) = - \langle E_\alpha, R(E_\beta,E_\gamma) \, E_\delta \rangle_T \, c^\beta c^\delta \, a^\gamma(r).
$$
Since
$$
-2\sum_{\alpha} \big( {a'_{\alpha}}^2+a_{\alpha}^2 \big) 
\leq
\frac{d}{dr} \Big( \sum_{\alpha}{a'_{\alpha}}^2+a_{\alpha}^2\Big)
\leq
2\sum_{\alpha} \big( {a'_{\alpha}}^2+a_{\alpha}^2 \big),
$$
we obtain $|V'(r)|_{T}\leq e^{r}$ and thus $|V(r)|_{T}\leq (e^r - 1)$.

By substituting this result into the above formulas, the
estimate can be improved again. Indeed, by computing and estimating the second-order derivative
$\frac{d}{dr}\sum_{\alpha}a_{\alpha}'a_{\alpha} $ as we did for the Jacobi
field estimate of Section~\ref{proofinject}, we can check that
$$
r - C(n) \, r^2 \leq \big( \sum|a_\alpha|^2(r) \big)^{1/2}\leq (e^r-1) \quad \text{ along the geodesic.} 
$$
Denote by $g_0, g_{T,0}$ the Lorentzian and the Riemannian metrics at the origin $0$
(which are nothing but the metrics at the point $p$),
and let $y^0,\ldots,y^n$ be Cartesian coordinates on $B_T(0,1)$,
with $\langle \frac{\del}{\del y^\alpha},\frac{\del}{\del y^\beta}\rangle_{g_0}(0) = \eta_{\alpha\beta}$.
Assuming that the radius under consideration is sufficiently small so that $(1 - C(n) \, |y|) < 1$
we conclude from the Jacobi field estimate that the exponential map is non-degenerate and that the metric
along the geodesic are comparable. In turn, since this is true for every radial geodesic, 
we can define the pull back of the metric to the tangent space and 
the conclusion hold in the whole ball $B_T(0,1)$, that is
\be
\label{metric7}
(1 - C(n) \, |y|) \, g_{T,0} \leq g_{T,y} \leq (1 + C(n) \, |y|) \, g_{T,0},
\qquad
y\in B_T(0,1).
\ee

By construction of the metric $g_T$ we have $\nabla_{g_T} - \nabla_g = \nabla_g T \ast T$ (schematically) and
$\nabla T(0) = 0$, and it is useful to control the covariant derivative too. To this end,
write the radial vector field as
$$
\frac{\del}{\del r}=\frac{y^{\alpha}}{r}\frac{\del}{\del y^{\alpha}},
\qquad
r :=  \big( \sum|y^{\alpha}|^2 \big)^{1/2},
$$
with $|\frac{\del}{\del r}|_{T}^2 \equiv 1$ (as stated already in \eqref{gamm}). 
Using that
$|\nabla T|_{T}^2 = \nabla_{\alpha}T^{\xi}\nabla_{\beta}T^{\eta} g_{T,\xi\eta} \, g_T^{\alpha\beta}$
and
computing the derivative of $|\nabla T|_{T}^2$ along radial geodesics, we find
$$
\frac{d}{dr}|\nabla T|_{T}^2
\leq C(n) \, |\nabla T|_{T}^3 + 2 \, \langle \nabla_{\frac{\del}{\del r}}\nabla T,\nabla T\rangle_{T}.
$$
By using that
$$
\nabla_{\frac{\del}{\del r}}T=0, \qquad
[\frac{\del}{\del r},\frac{\del}{\del y^{\alpha}}]
=-\frac{1}{r}\frac{\del}{\del y^{\alpha}}+\frac{y^\alpha}{r^2}\frac{\del}{\del r},
$$
we obtain
$$
\nabla_{\frac{\del}{\del r}}\nabla_{\frac{\del}{\del y^{\alpha}}}T^{\gamma}=-\frac{1}{r}\nabla_{\frac{\del}{\del
y^{\alpha}}}T^{\gamma}+R(\frac{\del}{\del r},\frac{\del}{\del y^{\alpha}})T^{\gamma},
$$
and therefore, thanks to the curvature assumption, 
$$
\frac{d}{dr}|\nabla T|_{T}^2\leq -\frac{2}{r}|\nabla T|_{T}^2
+ C(n) \, |\nabla T|_{T}^3 + C(n) \, |\nabla T|_{T}.
$$
This implies the following bound for the covariant derivative
\be
\label{covar7}
|\nabla T|_{T}(y)\leq C(n) \, |y|, \qquad |y|\leq 1/C(n),
\ee
which also provides a bound for the difference $\nabla_{g_T} - \nabla_g$.

\

\noindent{\sl Step 2. Estimate of the injectivity radius of $g$ on $B_T(0,c(n))$.}

Since the curvature on $B_T(0,1)$ is bounded and that $|\nabla_{g_T} - \nabla_g|^2_{T}\leq C(n)=1/c(n)$ on
the ball $B_T(0,c(n))$ we can follow
the argument in Section~\ref{proofinject} and bound from below the conjugate radius
for any point in the ball $B_T(0,3c(n)/4)$.

Next, given any point $y\in B_T(0,c(n)/2)$, let $\gamma_1$ and $\gamma_2$ be two geodesics
which meet at their end points and have ``short'' length with respect to the metric $g_T$ (or $g_{T,0}$).
By using the linear structure on $B_T(0,1)$ (a subset of the vector space $T_pM$) 
we can 
construct an homotopy of the loop $\gamma_1\cup\gamma_2^{-1}$
to the origin,
such that each curve have also ``short'' length for the metric $g_T$. By lifting the homotopy to the tangent space
$T_y B_T(0,1)$ and by relying on the conjugate radius bound, we reach a contradiction as
was done in Section~\ref{proofinject}.

In summary, there exists a universal constant $C(n)=1/c(n)$ depending only on the dimension such that the
injectivity radius at each point $y$ of $B_T(0,c(n))$ is
bounded from below by $4 c(n)$. Moreover, by the Jacobian field estimate again, we can prove the ball
$B_{T,p}(0,c(n))\subset T_pM$ defined by the Euclidean metric $g_{T,p}$ is covered by
$\exp_{y}(B_{T,y}(0,3c(n)))$, where $B_{T,y}(0,3c(n))\subset T_yT_pM$ is a ball of radius
$3c(n)$ defined by metric $g_{T,y}$, and any two points in
$B_{T,p}(0,c(n))$  can be connected by a $g$ geodesic which is totally contained in $B_{T,p}(0,2c(n))$. 
Further arguments are now required to arrive at the desired bound \eqref{toprove}.

\

\noindent{\sl Step 3. Riemannian metric $g_N$ induced on $B_T(0,2c(n))$. } Consider a geodesic $\gamma$ satisfying
$\gamma(0)=0$ and $\gamma'(0)=-T$, and let us set
$$
\gamma(c(n)/2)=:q, \qquad 
\tau := d_g(\cdot,q) - d_g(q,0). 
$$
Then, by following exactly the same arguments as in the main proof of Section~\ref{convexsection}, we construct a normal
coordinate system (of definite size) such that $g=-d\tau^2+g_{ij} \, dx^i dx^j$ and 
$g_N = d\tau^2+g_{ij}dx^{i}dx^j$,
and such that the corresponding reference Riemannian metric
satisfies the following properties:
\begin{itemize}
\item[(i)] $(1-c(n)) \, g_N \leq g_T \leq (1+c(n)) \, g_N$,
\item[(ii)] $g_N$ has bounded curvature ($\leq C(n)$) (see \eqref{timederiv}), and
\item[(iii)] for any fixed $y_0\in B_T(0,c(n))$ the distance function $d_{g_N}(y_0,)^2$ 
is strictly $g$-convex on the ball $B_T(0,2c(n))$ and, more precisely,
$$
(2+c(n)) \, g_N \geq \nabla^2_g d_\gt^2(y_0,\cdot) \geq (2-c(n)) \, g_N
\qquad \text{ on } B_T(0,2c(n))
$$
for any $y_0\in B_T(0,c(n))$.
\end{itemize}
The Hessian of the distance function (defined by the Riemannian metric $g_N$)
is naturally computed using the covariant derivatives defined by the Lorentzian metric $g$. 

\

\noindent{\sl Step 4. } Suppose that $p_1,\cdots,p_N$ are distinct
pre-images of $p$ in the ball $B_T(0,c(n))$. We claim
that any $p'\in \Bcal_T(p,c(n))$ has at least $N$ distinct pre-images in $B_T(0,1)$, 
and refer to this property as a ``lower semi-continuity" property.

Generalizing the terminology in \cite{CGT}, we use the
notation $a\mathop{\sim}\limits_{(g_{T,0},A)}b$ when two curves $a,
b$ defined on $M$ and with the same endpoints are homotopic
through a family of curves whose lift have $g_{T,0}$-lengths $\leq
A$. Relying on the lift and the linear structure, we see that, for
any curve $\xi$ starting from $p$ with (after lifting through $0$)
$g_{T,0}$-length $A\leq 1$, there exists a unique $g$-geodesic
$\gamma_\xi$ (with the same end points as $\xi$) defined on $M$
such that $\xi \mathop{\sim}\limits_{(g_{T,0},A)} \gamma_{\xi}$.
This fact establishes a one-to-one correspondence between the following three concepts: 
\par
(i) equivalence class of curves through $p$  with
$g_{T,0}$-lengths $\leq 3c(n)$, 
\par
(ii) radial geodesic
segments of $g_{T,0}$-lengths $\leq 3c(n)$, and 
\par 
(iii)
points in the ball $B_T(0,3c(n))\subset T_p M$.

Let $\sigma$ be a $g$-geodesic connecting $p$ to $p'$ in
$\Bcal_T(p,c(n))$. Observe that the images of the lines
$\overline{Op_i}$ by the exponential map,
$\sigma_i=\expb_p(\overline{Op_i})$, are distinct geodesic loops
through $p$. Denote by $\sigmat_i$ the lift of $\sigma_i\cup
\sigma$ through 0, and denote by $p_i'$ the end point of
$\sigmat_i$. Then, it is clear that all the points $p_i'$
($i=1,\cdots,N$) are the pre-images of $p'$ in $B_T(0,1/2)$. We
claim that they are distinct.

Indeed, assuming that $p_i'=p_j'$ for some $i\neq j$, we would
find $\sigma\cup \sigma_{i}
\mathop{\sim}\limits_{g_{T,0},2c(n)}\sigma\cup \sigma_j$, which
gives
$$
\sigma_i\mathop{\sim}\limits_{g_{T,0}, 3c(n)} \sigma^{-1}\cup
\sigma\cup\sigma_{i}\mathop{\sim}\limits_{g_{T,0},3c(n)}
\sigma^{-1}\cup\sigma\cup\sigma_j
\mathop{\sim}\limits_{g_{T,0},3c(n)}\sigma_j.
$$
This would imply
$\sigma_i\mathop{\sim}\limits_{g_{T,0},3c(n)}\sigma_j$ and,
therefore, $p_i=p_j$, which is a contradiction. In short, this
argument shows that the ``cancellation law'' holds for the
homotopy class of ``not too long'' curves.

\

\noindent{\sl Step 5. } Suppose that there exist two distinct
$g$-geodesics $\gamma_1:[0,l_1]\rightarrow M$ and
$\gamma_2:[0,l_2]\rightarrow M$ satisfying
$$
\gamma_1(0)=\gamma_2(0)=p,
\qquad 
|\gamma'(0)|_{T}^2=|\gamma'(0)|_{T}^2=1, 
$$
and meeting at
their endpoints, that is: $\gamma_1(l_1)=\gamma_{2}(l_2)$. Then,
let $l := l_1+l_2$ and $\gamma := \gamma_2^{-1}\cup
\gamma_1:[0,l]\rightarrow M$. Our aim is to prove that
$$
l \geq c(n) \, \vol_g\big( \Bcal_T(p,c(n)) \big),
$$
which will give us the desired injectivity radius.

From the loop $\gamma$ we define a map
$\pi_{\gamma}:B_T(0,c(n))\rightarrow B_T(0,2c(n))$
as follows: for any $y\in B_T(0,c(n))$, the point
$\pi_{\gamma}(y)$ is the end point of the lift $
\expb_p(\overline{Oy})\cup\gamma$ (through the origin). If one
would have $\pi_\gamma(y)=y$ then by the cancellation law
established in Step~4, we would have $\gamma
\mathop{\sim}\limits_{g_{T,0},2c(n)}0$, which is a contradiction.
So, the map $\pi_\gamma$ has no fixed point.

Without loss of generality, we assume that $l\leq c(n)^5$.
Let $N=[c(n)^3 /l]$ be the largest integer less than $c(n)^3 /l$, 
and let us use the notation $2 \gamma = \gamma \circ \gamma$, etc. 

\

\noindent{\bf Claim. } The classes 
$[\gamma],[2\gamma], \cdots,[N\gamma]$ are distinct homotopy
classes for the relation $\mathop{\sim}\limits_{g_{T,0},c(n)^2}$.

\

If this were not true, then by the cancellation law we would have
$[j\gamma] \mathop{\sim}\limits_{g_{T,0},c(n)^2}0$ for some $1\leq j\leq N$.
We already know that all $\pi_{\gamma}^{i}$ is defined from
$B_T(0,c(n)^2)$ to $B_T(0,c(n))$ for $i\leq j$. Since for any $y\in B_T(0,c(n)^2)$ we have
$$
\expb_p(\overline{Oy})\cup j\gamma
\mathop{\sim}\limits_{g_{T,0},c(n)}\expb_p(\overline{Oy}),
$$
which implies that $\pi_{\gamma}^j=id$. We use here the notation $\pi^2_\gamma := \pi_\gamma \circ \pi_\gamma$, etc. 

Then, we define a function $u:B_T(0,c(n))\rightarrow \R$ by
$$
u(y)=d_{\gt}^2(0,y)+ d_{\gt}^2(0,\pi_\gamma y) + \cdots+ d_{\gt}^2(0,\pi_\gamma^{j-1} y).
$$
Since $\pi^j_{\gamma}=id,$ it is easy to see $u(\pi_\gamma
y)=u(y)$ for any $y\in B_T(0,c(n))$. That is to say, $u$
is $\pi_\gamma$-invariant. By Step~3,
$u$ is strictly $g$-geodesically convex on $B_T(0,c(n))$. More precisely, since
for any $g$-geodesic $\xi:[0,s_0]\rightarrow B_T(0,c(n)),$ $\pi_{\gamma}^{i}\xi$ are still
$g$-geodesics in $B_T(0,c(n))$, and
\begin{equation*}
\begin{split}
\frac{d^2}{ds^2}u(\xi(s))&=\nabla^2d^2_{\gt}(0,\cdot)(\xi'(s),\xi'(s))
+\cdots+\nabla^2d^2_{\gt}(0,\cdot)
\big( {d\pi^{j-1}_{\gamma}}_{\xi(s)}(\xi'(s)),{d\pi^{j-1}_{\gamma}}_{\xi(s)}(\xi'(s)) \big)
\\
&\geq \gt(\xi'(s),\xi'(s)) > 0.
\end{split}
\end{equation*}
Observe that
$$
u\mid_{ B_T(0,c(n))^{c}}\geq
j(1-c(n))^2(c(n) - \frac{2lc(n)^3}{ l})^2 \geq \frac{jc(n)^2}{2},
$$
and
$$
u(0)\leq j(jl)^2\leq j \, c(n)^5 < \frac{j c(n)^2}{2},
$$
so the minimum of function $u$ over $\overline{B_T(0,c(n))}$ is
only achieved at at an interior point, say $y_0\in B_T(0,c(n))$.
Then by $\pi_\gamma$ invariance of $u$, we have
$u(\pi_{\gamma}y_0)=u(y_0) < j c(n)^2 / 2$, and this implies
$\pi_\gamma(y_0)\in B_T(0,c(n))$. By the injectivity radius estimate
at $y_0\in (T_pM,g)$, there exists a $g$-geodesic connecting $y_0$ 
to $\pi_\gamma(y_0)$, which is contained in $B_{T,p}(0,2c(n))$. By using the strong $g$-geodesic
convexity of $u$, we conclude that $\pi_{\gamma}y_0=y_0$. This
contradicts the fact that $\pi_\gamma$ has no fixed point, and the claim is proved.

\

\noindent{\sl Step 6. } The pull back of the volume element of
$g$ is the same as the one of $g_T$. By combining this observation with our results in Steps~4 and 5 we find
$$
\vol_{g_T}(B_T(0,1))\geq \frac{c(n)^3}{l} \, \vol_g(\Bcal_T(p,c(n))),
$$
which gives
$$
l \geq c(n) \, \frac{\vol_g(\Bcal_T(p,c(n)))}{\vol_{g_T}(B_T(0,1))}
\geq 
c(n) \, \vol_g(\Bcal_T(p,c(n))).
$$
The proof of Theorem~\ref{nofoliation} is completed.
\end{proof}


\section{Volume comparison for future or past cones}
\label{conclude1}

In Riemannian geometry, under a Ricci curvature lower bound,
Bishop-Gromov's volume comparison theorem allows one to compare the volume of small
and large balls in a sharp and qualitative manner. Let us return to Step~2 of Section~\ref{convexsection},
where we introduced the index form associated with the synchronous coordinate system on time-like geodesics.
By noticing that the index form is symmetric and that Jacobi fields minimize the index form (in some sense), 
we can extend the method of proof of the index comparison theorem. However, in a general Lorentzian manifold,
since the index form we needed (without imposing a restriction on the geodesics) is non-symmetric,
we need to adapt the method of the index comparison theorem, as follows.

\begin{theorem}[Volume comparison theorem for cones]
\label{volvol}
Let $(M,g)$ be a globally hyperbolic, Lorentzian ${(n+1)}$-manifold.
Fix $p\in M$ and a vector $T\in T_pM$ with $g_p(T,T)=-1$, and suppose that the
exponential map $\expb_p$ is defined on the ball $B_T(0,r_0) \subset T_pM$ 
(determined by the reference inner product $g_T$ at $p$). 
Suppose also that the Ricci curvature satisfies on $\Bcal_T(p,r_0)$
$$
\Ric_g(V,V) \geq -n \, K_2 \, | \, |V|_g^2 | \quad \text{ for all time-like vector fields } V. 
$$
Then for any $0<r<s<r_0$ the inequality
$$
\frac{\vol_g(\FC(p,r))}{\vol_g(\FC(p,s)))}\geq \frac{\vol_{K_2}(B(r))}{\vol_{K_2}(B(s)))}
$$
holds, with $\FC(p,r) := \expb_p(FC(p,r))$ and 
$$
FC(p,r) := \big\{ 0<|V|_{g_{T,0}}<r_0, \, |V|_g^2<0, \, \langle T,V\rangle_{g_{T,0}} < 0 \big\}
$$
and $\vol_{K_2}(B(r))$ is the volume of the ball with radius $r$, $B(r)$ (analogous to $\Bcal_T(p,r)\subset M$), 
in the simply-connected Lorentzian  ${(n+1)}$-manifold with constant curvature $K_2$
(that is, $R_{\alpha\beta\gamma\delta}=-K_2 \, (g_{\alpha\gamma
g_{\beta\delta}}-g_{\alpha\delta}g_{\beta\gamma}))$.

More generally, if $\Sigma$ is a subset in unit sphere $S^n$ such that
$|V|^2_{g}<0, \langle T, V\rangle_{g}<0$ for all $V \in \Sigma$,
then the inequality 
$$
\frac{\vol_{g}(\FC_{\Sigma}(p,r))}{\vol_{g}(\FC_{\Sigma}(p,s)))}\geq
\frac{\vol_{K_2}(B(r))}{\vol_{K_2}(B(s)))}
$$
holds with $\FC_{\Sigma}(p,r) : = \expb_p(FC_\Sigma(p,r))$ and 
$$
FC_\Sigma(p,r) := \big\{ V \in FC(p,r) \, / \, \frac{V}{|V|_{g_{T}}}\in \Sigma \big \}. 
$$
\end{theorem}

This result will be used shortly to control the injectivity radius of null cones, but is also 
of independent interest. For definiteness we state the result for future cones. 

\begin{proof}
Let $\gamma:[0,s_0] \to M$ be a future-oriented time-like geodesic satisfying $\gamma(0) = p$
and $|\gamma'(0)|_{g_{T}}=-1$.
We are going to use the standard technique to compute the rate of change of the volume element along
$\gamma$. Given $s_1\in (0,s_0)$ assume that any point in $(0,s_1]$ is
neither a conjugate point nor a cut point with respect to $p$.
Let $v_0=\gamma'(s_1),v_1,v_2,\cdots,v_n$ be an orthonormal basis at
$\gamma(s_1)$ with respect to $g_{\gamma(s_1)}$.
Let also $V_{\alpha}$ be the Jacobi field defined on $[0,s_1]$  and 
satisfying $V_{\alpha}(0)=0$
and $V_{\alpha}(s_1)=v_\alpha$.
Clearly $V_0=(s/s_1) \, \gamma'$, and the vectors $V_i$ and $\nabla_{\gamma'}V_i$ are orthogonal to $\gamma'$ for all $i\geq1$.

Consider the Jacobian of the exponential map $\varphi(s) := J(d\expb_{\gamma(s)})$, which is given by 
$$
\varphi(s)^2=\frac{|\gamma'(s)\wedge V_1(s)\wedge \cdots \wedge
V_n(s)|_{g}^2}{s^{2n} \, |\gamma'(0)\wedge V_1'(0)\wedge\cdots V_n'(0)|_{g}^2}. 
$$
Denote by $\varphi_{K_2}(s)$  the
corresponding quantity in the simply connected Lorentzian ${(n+1)}$-manifold with
constant curvature $-K_2$. Define the index form
$$
I_{s_1}(X,Y) := 
\int_{0}^{s_1} \Big( \langle\nabla_{\gamma'}X,\nabla_{\gamma'}Y\rangle_{g}-
\Riem_g(\gamma',X,\gamma',Y) \Big) \, ds,
$$
where $X,Y$ are vector fields along $\gamma$ and
$\Riem_g(\gamma',X,\gamma',Y) :=-\langle \Riem_g(\gamma',X)\gamma',Y\rangle_{g}$. Observe that $I_{s_1}$ is
symmetric in $X,Y$. It is easy to see 
\begin{equation*}
\begin{split}
\frac{d}{ds}\Big|_{s=s_1}\log \varphi^2&=\sum_{i}\langle
V_i'(s_1),V_i(s_1)\rangle_{g}-\frac{2n}{s_1}\\
&=\sum_{i}I_{s_1}(V_i,V_i)-\frac{2n}{s_1}.
\end{split}
\end{equation*}

Let $E_i(s)$ be the parallel transport of $v_i$ along $\gamma$. Since there are no conjugate points 
along $\gamma,$ the Jacobi field minimizes the index form among all vector fields with fixed
boundary values. This is the same as in Riemannian geometry. The 
reason is that the length of time-like geodesic without conjugate 
points is locally maximizing among all nearby time-like curves with the same end points.
Let $\Vt_i(s)=\frac{\sinh s}{\sinh s_1}E_i(s),$ then $ I_{s_1}(V_i,V_i)\leq
I_{s_1}(\Vt_i,\Vt_i) $ and
\begin{equation*}
\begin{split}
\frac{d}{ds}|_{s=s_1}\log \frac{\varphi^2}{(\varphi_{K_2})^2}
& \leq - \sum_{i} \int_{0}^{s_1}\frac{(\sinh  s)^2}{(\sinh s_1)^2}(\Riem_g(\gamma',E_i,\gamma'E_i)-K_2)
\\
&=-\int_{0}^{s_1}\frac{(\sinh s)^2}{(\sinh s_1)^2} (\Ric_g(\gamma',\gamma') - n \, K_2) \, ds \, \leq 0.
\end{split}
\end{equation*}

The following is a simple but very important observation due to Gromov, which 
we now extend to a globally hyperbolic Lorentzian manifold.
Let $A$ be the star-shaped domain (with respect to $0$) in $T_pM$, such that
$\expb_p: A \cap B_T(0,r_0)$ is a diffeomorphism on its image and the
image of $\del A \cap B_T(0,r_0)$ is set of cut locus (in $\Bcal_T(p,r_0)$).
Let $\chi_A$ be the characteristic function of $A$. Since $\varphi(s) / \varphi_{K_2}(s)$
is decreasing in $s$ we see that $\chi_A \varphi / \varphi_{K_2}$ is also decreasing in $s$.
Now, we get two functions on $B_T(0,r_0)$, whose quotient is decreasing along radial geodesics.
Observe that $M$ is globally hyperbolic, so any point in $\FC(p,r_0)$ is connected to $p$ by
a maximizing time-like geodesic. This also implies that the integration of $\chi_A \varphi$ over $B_T(0,s)$
gives the volume $\vol_{g}(\Bcal_T(p,s))$.
Then, by integrating $\chi_A \varphi$ and $\varphi_{K_2}$ over $B_T(0,s)$ and after
a simple calculation we deduce that $\vol_g(\FC(p,s)) / \vol_{K_2}(B(s))$ is decreasing in $s$.
The case of the ratio
$\vol_g(\FC_{\Sigma}(p,s)) / \vol_{K_2}(B(s))$ is similar.
The proof of the theorem is completed.
\end{proof}

We are now in a position to prove :

\begin{corollary}[Injectivity radius based on the volume of a future cone] 
\label{coro47}
Let $M$ be a manifold satisfying the assumptions in Theorem~\ref{nofoliation} and assumed to be globally hyperbolic,
and let $T\in T_pM$ be a reference vector.
Let $\Sigma$ be a subset in the unit sphere $S^n$ included in the future cone $N_p^+$.
If $\vol_g(\FC_{\Sigma}(p,r_0))\geq v_0>0$, then the inequality 
$$
{\Inj_g(M,p,T) \over r_0} \geq c(\Sigma) \, \frac{v_0}{r_0^{n+1}}
$$
holds, where
$\FC_{\Sigma}(p,r_0) := \expb_p( FC_p(r_0))$ with 
$$
FC_p(r_0) := 
\big\{ 0 < |V|_{T}<r_0, \,
\langle T,V\rangle_{T}<0, \, |V|_{g}^2<0, \, \frac{V}{|V|_{T}}\in
\Sigma \big\},
$$
and the constant $c(\Sigma)$ depends only on the distance (measured by $T$) of $\Sigma$ to the null cone.
\end{corollary}

\begin{proof} First, we recall there is a constant $C(\Sigma)$ depending only on the distance of $\Sigma$ to the null
cone,  such that $\Ric(\gamma',\gamma')\geq -C(\Sigma)|\gamma'|_{g}^2$ for any time-like geodesic $\gamma$
with $\gamma'(0)\in \Sigma$. By the volume comparison theorem for future cone established 
in Theorem~\ref{volvol} we have
$$
\frac{\vol_g(\FC_{\Sigma}(p,c(n) r_0))}{\vol_g(\FC_{\Sigma}(p,r_0))}\geq
C(\Sigma),
$$
and combining this result with Theorem~\ref{nofoliation}, the corollary follows.
\end{proof}


\section{Final remarks}
\label{conclude2}

\subsection*{Regularity of Lorentzian metrics}

Following the strategy proposed in the present paper, we now ``transfer'' to the Lorentzian metric the
regularity available on a reference Riemannian metric. Clearly, the regularity obtained in this manner depends on the way
the reference Riemannian metric is constructed. The interest of our approach below 
is to provide a simple derivation: using harmonic-like coordinates for the Riemannian metric
we see immediately that the Lorentzian metric has uniformly bounded first-order derivatives.
For the optimal regularity achievable with Lorentzian metrics we refer to Anderson \cite{Anderson2}.

\begin{proposition}[Regularity in harmonic-like coordinates]
Under the assumptions and notation of Theorem~\ref{nofoliation}, define 
$$
r_1 := c(n) \, \frac{\vol_g(\Bcal_T(p, c(n) \, r_0))} {r_0^{n+1}} \, r_0,
$$
where $c(n)$ is the constant determined in this theorem. 
Then for any $\eps>0$ there exist a constant $c_1(n,\eps)$ with $\lim_{\eps\rightarrow0} c_1(n,\eps) = 0$
and a coordinate system $(x^\alpha)$ satisfying $x^{\alpha}(p)=0$
and defined for all $(x^0)^2 + (x^1)^2 + \ldots + (x^n)^2 < (1-\eps)^2 \, r_1^2$, such that
\be
\label{result}
\aligned
& |g_{\alpha\beta}-\eta_{\alpha\beta}|\leq c_1(n,\eps),
\\ 
& r_1 \, |\del g_{\alpha\beta}|\leq c_1(n,\eps), 
\endaligned
\ee
where $\eta_{\alpha\beta}$ is the Minkowski metric in these coordinates. 
\end{proposition}

\begin{proof}
By scaling we may assume $r_1=1$. By Step~1 in the proof of Theorem~\ref{nofoliation},
we know that the Riemannian metric $g_T$ is equivalent to the Riemannian metric $g_{T,0}$ on
$B_T(0,4 c(n))$. By considering a lift and using again the results in Step~1 this implies 
$$
\Bcal_T(p,c(n))\subset \Bcal_T(q,3 c(n))
\quad
q \in \Bcal_T(p, c(n)).
$$
Applying the same argument as in Theorem~\ref{nofoliation}, we deduce that the injectivity radius
of any point in $\Bcal_T(p,c(n))$ is bounded from below by $c(n)$.
As in Step~3 in the proof of Theorem~\ref{nofoliation} (or in Step~2 of Section~\ref{convexsection}),
we see that there exists a synchronous coordinate system $(y^\alpha) = (\tau, y^j)$ 
of definite size around $p$
such that the metrics $g = - d\tau^2 + g_{ij} \, dy^i dy^j$ and
$g_N = d\tau^2 + g_{ij} \, dy^i dy^j$ (the Riemannian metric constructed therein)
satisfy the following properties on the geodesic ball $\Bcal_T(p,c(n))$:
\begin{itemize}
\item[(a)] $(1-c(n)) \, g_N \leq g_T \leq (1 + c(n)) \, g_N$,
\item[(b)] $g_N$ has bounded curvature ($\leq 1/c(n)$),
\item[(c)] $|\tau| + {1 \over |\tau|} + |\nabla^2 \tau|_N \leq 1/c(n)$.
\end{itemize}
(In particular, this implies $|\nabla_{g_N} g|_N < 1/c(n)$.)  Since the volume
$\vol_{g_N}(\Bcal_T(p,c(n)))$ is bounded from below, it follows from \cite{CGT}
that the injectivity radius of $g_N$ at $p$ is bounded from below by $c(n)$.
By the theorem in \cite{JostKarcher} on the existence of harmonic coordinates,
for any small $\eps>0$ there exists an harmonic coordinate system $(x^\alpha)$ 
with respect
to the Riemannian metric $g_N$ such that $\sum_{\alpha}|x^{\alpha}|^2 < (1- \eps)^2$  and 
for every $0<\gamma<1$
$$
|g_{N,\alpha\beta} - \delta_{\alpha\beta}| < c_1(n,\eps),
\qquad
    |\del g_N| < 1/c(n),
    \qquad
    |\del g_N|_{\Cbf^{\gamma}} < 1/c(n,\eps,\gamma).
$$
In the construction of harmonic coordinates, 
we may also assume that
$|\frac{\del}{\del y_0}-\frac{\del}{\del \tau}|_{g_{T,0}} < c_1(n,\eps)$.

Since $|\nabla_{g_N} g|_N < 1/c(n)$ and that, in these coordinates, $|\nabla_{g_N}|\leq 1/c(n)$,
we have $|\del g|< 1/c(n)$. Finally, to estimate the metric we write
$|g_{\alpha\beta} - \eta_{\alpha\beta}|_p < c_1(n,\eps)$ and
 $|\del g| < 1/c(n)$ and we conclude that
$|g_{\alpha\beta} - \eta_{\alpha\beta}| < {1 \over C(n)} \eps + c_1(n,\eps)$.
The proof is completed.
\end{proof}


\subsection*{Pseudo-Riemannian manifolds}

Finally, we would like discuss pseudo-Riemannian manifolds $(M,g)$ (also referred to as semi-Riemannian manifolds).
Consider a differentiable manifold $M$ endowed with a symmetric, non-degenerate covariant $2$-tensor $g$.
We assume that the signature of $g$ is $(n_1,n_2)$, that is, $n_1$ negative signs and $n_2$ positive signs.
Riemannian and Lorentzian manifolds are special cases of pseudo-Riemannian manifolds.
Fix $p\in M$ and an orthonormal family $T$ consisting of $n_1$ 
vectors $E_1,E_2,\cdots,E_{n_1}\in T_pM$ such that 
$\langle E_i, E_j\rangle_{g}=-\delta_{ij}$. Based on this family, we can define a reference inner product
$g_T$ on $T_pM$ by generalizing our construction in the Lorentzian case, and by 
using this inner product we can then define the ball $B_T(0,s)\subset T_pM$. 
By parallel translating $E_1,E_2,\cdots,E_{n_1}$ along radial geodesics from the origin in $T_pM$, 
we obtain vector fields $E_1,E_2,\cdots,E_{n_1}$ defined in the tangent space 
(or multi-valued vector fields on the manifold). 
This also induces a (multi-valued) Riemannian metric $g_T$ as was explained before.

The following corollary immediately follows by repeating the proof of Theorem~\ref{nofoliation}.
We note that the curvature covariant derivative bound imposed below is probably superfluous and 
could probably be removed by introducing a foliation based on certain synchronous-type coordinates, 
as we did in Section~\ref{convexsection} for Lorentzian manifolds. On the other hand, to the best of our knowledge
this is the first injectivity radius estimate for pseudo-Riemannian manifolds.

\begin{corollary}[Injectivity radius of pseudo-Riemannian manifolds]
\label{coro}
Let $(M,g)$ be a 
differentiable pseudo-Riemannian $n$-manifold with signature $(n_1,n_2)$,
and let $p\in M$ and $T=(E_1,\cdots,E_{n_1})$ be a family of vectors in $T_pM$ 
satisfying $g(E_i,E_j) = -\delta_{ij}$.
Suppose that the exponential map $\expb_p$ is defined on $B_T(0,r_0)\subset T_pM$ and that
$$
| \Riem_g |_{T} \leq r_0^{-2},
\quad  |\nabla \Riem_g|_{T}\leq r_0^{-3} \qquad \text{ on } B_T(0,r_0).
$$
Then, there exists a positive constant $c(n)$ such that
$$
{\Inj_g(M,p,T) \over r_0}
\geq
c(n) \, \frac{\vol_g(\Bcal_T(p,c(n) \, r_0))}{r_0^n},
$$
where $\Bcal_T(p,r)=\expb_p(B_T(0,r))$ is the geodesic ball at $p$ with radius $r$.
\end{corollary}

\begin{proof}
Without loss of generality we assume $r_0=1$.  In local coordinate system $y^{\alpha}$,  let
$$
E_i =: E_i^\beta \frac{\del}{\del y^{\beta}},
\quad 
    E_{i\alpha}= E_{i}^{\beta}g_{\alpha\beta},  \qquad i=1, \ldots, n_1, 
$$
then
$g_{T,\alpha\beta}=g_{\alpha\beta}+2\sum_{i=1}^{n_1}E_{i\alpha}E_{i\beta}$.
By the same computations as in the proof of
Theorem~\ref{nofoliation} we obtain
$$
\aligned
& |\nabla E_i|_{T}\leq {1 \over c(n)},
\\
& |g_T - g_{T,0}| + |g-\eta| < c(n) \qquad \text{ on the ball } B_T(0,c(n)),
\endaligned
$$
where $\eta_{\alpha\beta} := \mp \delta_{\alpha\beta}$
(a minus sign for $\alpha\leq n_1$, and a plus sign for $\alpha>n_1$).
In view of the computations in \cite{Hamilton} (Theorem~4.11 and Corollary~4.12) we deduce
that 
$|\del g|< r/c(n)$, where $r^2 = (y^1)^2 + \cdots + (y^n)^2$. 
Since $d_{g_{T,0}}^2(y_0,y)=|y-y_0|^2$, we have for any point
$y_0\in B_T(0,c(n))$ 
\begin{align*}
\nabla^2_{\alpha\beta}d_{g_{T,0}}^2(y_0,\cdot) \geq \delta_{\alpha\beta}=g_{T,0} 
\qquad \text{ on the ball } B_T(0,c(n)).
\end{align*}
Since the metric $g_{T,0}$ plays the same role as $g_N$ (cf.~the proof of Theorem~\ref{nofoliation}),
all arguments can be carried out and this completes the proof of the corollary.
\end{proof}


\section*{Acknowledgements}
The first author (BLC) was partially supported by Sun Yat-Sen University via 
China-France-Russia collaboration grant (No.~34000-3275100),
the Ecole Normale Sup\'e\-rieure de Paris, the French Foreign Ministry,
and the Institut des Hautes Etudes Scientifiques (IHES, Bures-sur-Yvette).
The second author (PLF) was partially supported by the A.N.R. (Agence Nationale de la Recherche)
through the grant 06-2-134423
entitled {\sl ``Mathematical Methods in General Relativity''} (MATH-GR), and by the Centre National de la Recherche Scientifique (CNRS).


\end{document}